\documentclass[11pt,twoside]{amsart}
\usepackage{amssymb}
\usepackage{mathrsfs}
\def\flex#1{\mathrel{\mathop{\kern 0pt\hbox to 
10mm{\rightarrowfill}}\limits_{#1
\rightarrow \infty}}}

\def\Flex#1{\mathrel{\mathop{\kern 0pt\hbox to 
10mm{\rightarrowfill}}\limits_{#1
\rightarrow \infty}}}

\font\bb=msbm10 at 12pt

\def\P{\hbox{\bb P}}
\def\Q{\hbox{\bb Q}}

\def\1{1\hspace{-1.2mm}\mbox{{\normalsize I}}}

\newtheorem{Th}{Theorem}
\newtheorem{Def}[Th]{Definition}
\newtheorem{Lemme}[Th]{Lemma}
\newtheorem{Cor}[Th]{Corollary}
\newtheorem{Rmq}[Th]{Remark}
\newtheorem{Prop}[Th]{Proposition}
\newcommand{\R}{\mathbb{R}}
\newcommand{\E}{\mathbb{E}}
\newcommand{\F}{\mathcal{F}}

\title{A two-dimensional oblique extension of bessel processes}
\begin{document}
\maketitle
\centerline{Dominique L\'{e}pingle\footnote{Universit\'{e} d'Orl\'{e}ans, MAPMO-FDP, F-45067 Orl\'{e}ans 
(dominique.lepingle@univ-orleans.fr)}}

\begin{abstract}
We consider a Brownian motion forced to stay in the quadrant by an electrostatic oblique repulsion from the sides. We tackle the question of hitting the corner 
or an edge, and find product-form stationary measures under a certain condition, which is reminiscent of the skew-symmetry condition for a reflected 
Brownian motion.
\end{abstract} 

\section{Introduction}

In the present paper we study existence and properties of a new process with values in the nonnegative quadrant  $S=\R_+\times \R_+$ where 
$\R_+:=[0,\infty)$.
It may be seen as a two-dimensional extension of a usual Bessel process. It is a two-dimensional Brownian motion forced to stay in the quadrant by 
electrostatic repulsive forces, in the same way as in the one-dimensional case where a Brownian motion 
which is  prevented from becoming negative by an electrostatic drift becomes a Bessel process.
Note here and now that the corner  ${\bf 0}=(0,0)$ will play a crucial role 
and in some cases
it will be necessary to restrict the state space to the punctured nonnegative quadrant 
$S^{{\bf 0}}=S\setminus \{{\bf 0}\}$.

Let $(\Omega,\F, (\F_t)_{t\geq 0}, \P)$ be a complete probability space endowed with a filtration
$(\F_t)_{t\geq 0}$ satisfying the usual conditions. 
Let  $(B_t,C_t)$ be an adapted driftless  Brownian motion in the plane starting from ${\bf 0}$, 
with covariance matrix
\[
  \left( \begin{array}{cc}
        1 & \rho \\
        \rho & 1
        \end{array}
  \right) 
\]
and $\rho\in[-1,+1]$.

\begin{Def}
  \label{Th:defi}
 Let $\alpha,\beta,\gamma,\delta$ be four real constants with $\alpha>0,\delta>0$. We say that an
$(\F_t) $-adapted continuous process $(X,Y)$ with values in $S$ is an Oblique Two-dimensional Bessel Process (O2BP) if for any $t\geq 0$
\begin{equation}
 \label{eq:obli}
 \begin{array}{lllll}
 X_t & = & X_0+B_t+\alpha\int_0^t\frac{ds}{X_s} + \beta\int_0^t\frac{ds}{Y_s} & \geq & 0 \\
 Y_t & = & Y_0+C_t+\gamma\int_0^t\frac{ds}{X_s} + \delta\int_0^t\frac{ds}{Y_s} & \geq & 0
 \end{array}
\end{equation}
where $X_0$ and $Y_0$ are nonnegative ${\F}_0$-measurable random variables, and
\[
\begin{array}{lll}
  \int_0^t{\bf 1}_{\{X_s=0\}}ds =0 & \quad & \int_0^t{\bf 1}_{\{Y_s=0\}}ds =0\\
  \int_0^t{\bf 1}_{\{X_s>0\}}\frac{ds}{X_s}< \infty  & \quad & \int_0^t{\bf 1}_{\{Y_s>0\}}\frac{ds}{Y_s}< \infty \,.
\end{array}
\]
\end{Def}
 \vspace{0.3cm}

This stochastic differential system is very singular at the edges of the quadrant and the question of existence and uniqueness of a solution
is not simple.  The particular case  when $\beta=\gamma=0$ 
\begin{equation}
 \label{eq:normal}
   \begin{array}{lllll}
 U_t & = & X_0+B_t+\alpha\int_0^t\frac{ds}{U_s}  & \geq & 0 \\
 V_t & = & Y_0+C_t + \delta\int_0^t\frac{ds}{V_s} & \geq & 0 \,.
 \end{array}
\end{equation}
is already known:
the processes $U$ and $V$ are  Bessel processes.
Actually, $U$ is a Bessel process of dimension $2\alpha +1$, and the point $0$ is instantaneously reflecting for $U$ if 
$\alpha<\frac{1}{2}$ and polar if $\alpha\geq\frac{1}{2}$. If $\rho=0$, $U^2+V^2$ is the square of a Bessel process 
of dimension  $2\alpha+2\delta + 2$ \cite{SW73}, and so the corner ${\bf 0}$ is polar for $(U,V)$ in this case. Comparison between 
$X$ and $U$, $Y$ and $V$ will play a key role in the construction of the solution $(X,Y)$ and the 
study of its behavior close to the edges of the quadrant. The process $(U,V)$ is an example of  Brownian motion perturbed by a drift deriving 
from a convex potential. More generally, stochastic differential systems including such a singular drift have been studied in \cite{C95, P95, St95, R96}, where 
strong existence and uniqueness were obtained. They are examples of so-called multivalued stochastic differential equations, also called stochastic variational
inequalities in convex analysis. 

We will use these results to study  the solutions to (\ref{eq:obli}) in the oblique case where $\beta$ or $\gamma$ do not vanish. 
We obtain strong existence and uniqueness for a large 
set of parameters and initial conditions, but not for all possible values. 
In the proofs we naturally fall into the crucial question of hitting the corner, that is the non-smooth part of the boundary. Using McKean's argument 
 on the asymptotic behavior of continuous local martingales obtained by time change from the real driftless Brownian motion,  we are able to state several
sufficient conditions to prevent the processes from hitting the corner.  Our methods are not powerful enough to allow for necessary conditions.  However we shall
not restrict to processes avoiding the corner and, depending on the parameters, we will get existence and uniqueness (in a strong sense) sometimes in the
 whole quadrant, sometimes in the punctured quadrant, that is, the quadrant without a corner.

We will also obtain some partial results about the attainibility of the edges  of the quadrant. It is interesting to see whether
 the boundary behavior of one component may be modified by the interaction with the other component.

In the one-dimensional case, the so-called scale functions transform the solution of a stochastic differential equation into  local martingales and one may infer  some information on the boundary behavior of the process. 
This technique was very successful in the study of Bessel processes. Here we still obtain  functions of an O2BP which are local martingales or supermartingales for some values of the parameters and we derive some information on its  asymptotic 
behavior. 
 
The laws of Bessel processe with different parameters \cite{RY94} are mutually abolutely continuous when their paths do not reach the origin.  
Here we  obtain two partial results about the absolute continuity of the laws of 
O2BPs for some values of the parameters. 

Finally we follow the way explored in \cite{OO14}, where a drifted Brownian motion is endowed with another drift term that is
continuous and depends obliquely, via a regular potential function, on the position of the process relative to an orthant. Under an additional condition, 
which is called a skew-symmetry condition,  an invariant density was given in an explicit product form. We obtain the same result  for an O2BP (where
the Brownian motion is now drifted) under an anologous condition. This time, the invariant density is the product of two gamma densities, which is consistent 
with the one-dimensional case where the invariant measure of a drifted Bessel process is a gamma distribution.

\subsection{Comparison with obliquely reflected Brownian motion}

Studying O2BPs  makes appear a strong connection with the properties of a semimartingale reflecting  Brownian motion in the quadrant.
We briefly recall the definition. Let $(B_t,C_t)_{t\geq 0}$ be the Brownian motion in Definition \ref{Th:defi}.

\begin{Def}
\label{Th:refobli}
An semimartingale  reflecting Brownian motion (SRBM) in the quadrant is a continuous adapted process $(X_t,Y_t)$ which is a solution to the system
\begin{equation}
 \label{eq:broref}
   \begin{array}{lllll}
    X_t & = & X_0+B_t+L_t^1 + r_1 L_t^2 & \geq & 0 \\
   Y_t & = & Y_0 + C_t + r_2 L_t^1 + L_t^2 & \geq & 0
   \end{array}
\end{equation}
where  $L^1$ and $L^2$ are two continuous adapted nondecreasing processes with $L_0^1=L^2_0=0$ such that for any $t\geq 0$
\[
    \int_0^t {\bf 1}_{\{X_s>0\}}dL_s^1 = \int_0^t {\bf 1}_{\{Y_s>0\}} dL^2_s = 0.
\]
Here $X_0$ and $Y_0$ are nonnegative  ${\F}_0$-measurable random variables, $r_1$ and $r_2$ are real numbers.
\end{Def}

There was an extensive literary output on that topic in the eighties, with a more general domain : a wedge, an orthant or a convex polyhedron. 
We mention the works
 of Harrison, Reiman, Varadhan, Williams, Dai \cite{H78, HR81, VW85, W87, DW96} to cite a few of them. For a  more 
complete bibliography we refer to  \cite{W95}. 

We note a first analogy between their results and ours: we prove a necessary condition of existence of an O2BP is the existence of a convex combination of the 
directions of interaction
\begin{equation}
 \label{eq:repul}
 {\bf  r_x} = \left( \begin{array}{l} \alpha \\ \gamma \end{array}\right) \;\; {\rm and} \;\; 
{\bf  r_y} = \left(\begin{array}{l} \beta \\ \delta \end{array}\right)
\end{equation}
that points into the quadrant from the corner. In the reflection setting  it was proved in \cite{W85} that a necessary an sufficient condition for the existence 
of a SRBM  is the existence of a convex combination of the directions of reflection
\[
\left( \begin{array}{l} 1 \\ r_2 \end{array}\right) \;\; {\rm and} \;\; 
\left(\begin{array}{l} r_1 \\ 1 \end{array}\right) \;
\]
with the same property.

Another analogy is met in the delicate and important question of attainibility of the corner. The authors in \cite{VW85} have found an explicit harmonic function that provides a
full answer to the question of hitting the corner for a reflected Brownian motion in a wedge of angle $\xi \in (0,2\pi)$ with the identity matrix as covariance matrix.
We do not have any such convenient function. However, our second condition 
\begin{equation}
 \label{eq:evitebes}
  2\rho \leq \frac{\beta}{\delta}+\frac{\gamma}{\alpha}
\end{equation}
in Corollary \ref{Th:asym}
is reminiscent of the  necessary and sufficient condition in \cite{S15a} which writes 
\begin{equation}
 \label{eq:eviteref}
  2 \rho \leq r_1+r_2
\end{equation}
in the setting of Definition \ref{Th:refobli}. 
Restricting to $\xi\in(0,\pi)$, we may consider a linear transformation in the plane that changes the wedge in \cite{VW85} into the nonnegative quadrant and the 
initial Brownian motion into a Brownian motion with covariance coefficient $\rho= - \cos{\xi}$. Then the non-attainibility
condition in \cite{VW85} is transformed into the above condition (\ref{eq:eviteref}). 

For theoretical as well as practical reasons, a great deal of interest was taken in the question of recurrence of the Brownian motion with a constant drift vector and oblique reflection,  and in the computation of the invariant measure
\cite{W85a, WZ90, HH09, S15c}.
 Under the assumption that the directions of reflection satisfy a skew-symmetry condition,
it was proved that the invariant measure has exponential product form density \cite{HW87a, HW87, W87, W95}. 
This result has a practical interest because of the interpretation of the SRBM as approximation of the queue length processes for
networks of queues in heavy traffic \cite{H78, HW87a}. There is also a financial reason for studying SRBMs. 
Motivated by the so-called Atlas model of equity market
presented in \cite{F02}, 
some authors \cite{IK10, IP11, IKS11, IK12, IPS12, S15a} have recently studied competing  Brownian particles on the line with rank dependent 
local characteristics. In fact the gaps between adjacent particles are proportional to the components of a SRBM in an orthant. There is an
 invariant probability density with an explicit exponential product form when the volatility coefficients are constant \cite{PP08}. More
generally, this is still true if 
their squares depend on rank linearly since one may  infer from the statements in Section 2 of \cite{S15a} that the skew-symmetry condition 
is still satisfied. 
Following the way in \cite{HW87}, \cite{DW94} and \cite{HH09}, a more general question is 
 the recurrence or transience of an O2BP. 
 Answering this question does not seem to be an easy task. The method in \cite{HR93}, which provides a full answer for the obliquely reflected 
Brownian motion in the quadrant, appears to break down here.
 Mimicking the computation in \cite{OO14}, we  just calculate an invariant measure in product form under a  skew-symmetry condition which 
is the equality condition in the inequality (\ref{eq:evitebes}), whereas   the skew-symmetry condition  in the reflecting case is the equality condition 
in (\ref{eq:eviteref}).
 Now the terms of the product are  gamma
distributions with explicit parameters.

Another topic of interest in \cite{OO14} lies in its Remark 4.12: when a scale parameter  goes to zero, the exponentially reflected Brownian motion
 should converge to the 
obliquely reflected Brownian motion. It is typically a penalty method. This kind of approximation by a sequence of diffusions with regular drifts living on
the whole Euclidean space has been used for instance in  \cite{LS84,S15b} for oblique reflection in domains with smooth boundary
and in \cite{M83,Sl13} for normal reflection in convex domains. With regard to our framework, it can be shown that a sequence of real Bessel processes 
with dimension decreasing to 1 pathwise decreases to a  reflected real Brownian motion. It could be interesting to consider 
a sequence of oblique two-dimensional Bessel 
processes with interaction matrices
\[
  n^{-1} \left(\begin{array}{cc}
     1 & r_1 \\
     r_2 & 1 
     \end{array}
     \right)
\]
and to let $n$ go to infinity.  It should converge to the SRBM in Definition \ref{Th:refobli}. 
This time, this would be an interior approximation, as it was done in \cite{WZ90} and \cite{PW94} in the particular setting of stationary Markov processes 
associated with Dirichlet forms. We leave aside this point for further investigation.

\subsection{Organization of the paper}
The rest of the paper is organized as follows. In Section 2 we recall some trajectorial properties of usual Bessel processes. In Section 3 we state and
prove  three main lemmas of repeated use in the sequel. We also recall 
an existence and uniqueness result for a special case of multivalued stochastic differential equation that will be useful in our construction of an OB2P.
  Section 4 is devoted to the proof of sufficient conditions to avoid the corner of the quadrant.
 The main theorems of existence and uniqueness of an O2BP are given in Section 5. In Section 6 we discuss the question of hitting the edges of the quadrant. In Section 7 we consider two particular cases where there exist simple functions of O2BPs that are local martingales and obtain some information on the
 asymptotic behavior of the paths. We leave the trajectorial point of view in Section 8 to tackle questions of absolute continuity of the law of an O2BP with respect to the product of laws of real Bessel processes. The final Section 9  introduces a skew-symmetry condition that allows us to obtain existence of a stationary probability
in form of the product of two gamma distributions.

\vspace{0.5cm} 

\section{Some properties of Bessel processes}

A Bessel process of dimension $d>1$, starting at $r\geq 0$, is the unique solution to the stochastic differential equation
\begin{equation}
 \label{eq:bessel}
  R_t = r + W_t + \frac{d-1}{2} \int_0^t \frac{1}{R_s}ds .
\end{equation}
where $W$ is a standard driftless  real-valued Brownian motion starting at $0$.

We know that:
\begin{itemize}
\item
{\bf $d\geq 2$} : the  point $0$ is a  polar set (\cite{RY94}, Proposition V.2.7);
\item
{\bf $d=2$} : $\limsup_{t\rightarrow \infty} R_t= +\infty,\; \liminf_{t\rightarrow \infty} R_t = 0$ (\cite{RY94}, Theorem V.2.8);
\item
{\bf $1<d<2$} : the point $0$ is instantaneously reflecting (\cite{RY94}, Proposition XI.1.5).
\end{itemize}

More precisely (\cite{JYC09}, p.337 and p.339),
\begin{itemize}
\item
{\bf Dimension $d>2$}: 
\[
 \begin{array}{l}
   \P(R_t>0, \forall t>0)=1, \\
   \P(R_t\rightarrow \infty, t\rightarrow \infty)=1, \\
   R^{2-d} \;\rm{ is\; a\; local\; martingale}.
 \end{array}
\]
\item
   {\bf Dimension $d=2$}:
  \[
  \begin{array}{l}
   \P(R_t>0, \forall t>0)=1, \\
   \P(\sup_{t>a} R_t=\infty,\inf_{t>a} R_t=0)=1 \;\; {\rm  for \; any }\; a>0,\\
   \ln R\; \rm{ is\; a\; local\; martingale}.
  \end{array}
 \]
\item
  {\bf Dimension $1<d<2$}:
  \[
    \begin{array}{l}
     P(R_t>0, \forall t>a)=0 \;\;{ \rm for\; any }\; a>0.
    \end{array}
   \]
 \end{itemize}

\begin{Cor}
 \label{Th:besse}
When $d\geq 2$, for any $a>0$,
\[
   \int_a^{\infty} \frac{ds}{R_s^2} = \infty \;\, a.s.
\]
\end{Cor}
{\bf Proof}.
From It\^{o}'s formula
\[
  \begin{array}{lll}
    \ln{R_t} & = & \ln{R_a} + \int_a^t\frac{dW_s}{R_s} + \frac{d-2}{2} \int_a^t \frac{ds}{R_s^2}.
  \end{array}
\]
From (\cite{RY94}, Proposition IV.1.26) we know that on the set $\{\int_a^{\infty}\frac{ds}{R_s^2}<\infty\}$, the continuous 
local martingale $\int_a^t\frac{dW_s}{R_s}$ converges as $t\rightarrow \infty$. But we have seen that with unit probability  
$\ln{R_t}$ does not converge.
$\hfill \blacksquare$

\vspace{0.3cm}

We shall also use an absolute continuity result (\cite{RY94}, Exercise XI.1.22 or \cite{JYC09}, Proposition 6.1.5.1). On the canonical space $\Omega =
C([0,\infty),[0,\infty))$, we denote by $R$
 the canonical map $R_t(\omega)=\omega (t)$, by ${\mathcal R}_t=\sigma(R_s, s\leq t)$ the canonical filtration and by 
$\P^d_r$ the law of the Bessel process of dimension $d\geq 2$ starting at $r>0$. Then, 
\begin{equation}
 \label{eq:absocon}
\P^d_r\mid_{{\mathcal R}_t}= \left(\frac{R_t}{r}\right )^{\frac{d-2}{2}} \exp{\left (-\frac{(d-2)^2}{8}\int_0^t \frac{ds}{R_s^2}\right )} \, . \, 
\P^2_r\mid_{{\mathcal R}_t} .
\end{equation}

\vspace{0.5cm}

\section{Four useful tools}

 The following simple lemma is a comparison lemma. It will play an ubiquitous role in our proofs.
\begin{Lemme}
For $T>0$, $\alpha>0$, let $x_1$ and $x_2$ be nonnegative continuous solutions on $[0,T]$ to the equations
\[
  \label{Th:comp}  
\begin{array}{lll}
   x_1(t) & = & v_1(t) + \alpha\int_0^t\frac{ds}{x_1(s)} \\
   x_2(t) & = & v_2(t) + \alpha\int_0^t\frac{ds}{x_2(s)}
  \end{array}
\]
where  $v_1$, $v_2$ are continuous functions such that $0\leq v_1(0)\leq v_2(0)$, and $v_2-v_1$ is 
nondecreasing.
Then $x_1(t)\leq x_2(t)$ on $[0,T]$.
\end{Lemme}
{\bf Proof}.
Assume there exists $t\in (0,T]$ such that $x_2(t)<x_1(t)$. Set
  \[
   \tau:=\max \{s\leq t : x_1(s)\leq x_2(s)  \} \,.
 \]
Then,
\[
  \begin{array}{lll}
  x_2(t)-x_1(t) &= & x_2(\tau)-x_1(\tau)+(v_2(t)-v_1(t))-(v_2(\tau)-v_1(\tau)) +
  \alpha  \int_{\tau}^t(\frac{1}{x_2(s)}-\frac{1}{x_1(s)})ds \\
   &\geq & 0 \,,
  \end{array}
\]
a contradiction. $\hfill  \blacksquare$

\vspace{0.4cm}

The following elementary lemma will also be repeatedly used.

\begin{Lemme}
\label{Th:seconddegre}
Let $Q(x,y)=ax^2+bxy+cy^2$ be a  second degree homogeneous polynomial. Then $Q$ is nonnegative on the whole $S$ 
if and only if $a\geq 0$, $c\geq 0$ and $b\geq -2\sqrt{ac}$.
\end{Lemme}
{\bf Proof}. \\
Taking $x=0$ (resp. $y=0$) we see that $a$ (resp. $c$) must be nonnegative. Then we rewrite
\[
  Q(x,y) = (\sqrt{a}x-\sqrt{c}y)^2 + (b+2\sqrt{ac})xy .
\]
So $Q$ is nonnegative for all $x\geq 0$ and $y\geq 0$ iff $b\geq -2\sqrt{ac}$. $\hfill  \blacksquare$

\vspace{0.4cm}

Another main tool will be the following convergence result, whose argument goes back to McKean (\cite{MK69}, p.31 and p.47). The statement and the proof below
are borrowed from (\cite{RY94}, Theorem V.1.7 and Proposition V.1.8). The only change lies in the introduction of a stopping time $\tau$ up to which the local martingale $M$ is now defined.
\begin{Lemme}
\label{Th:kean}
Let $M$ be a continuous local martingale defined on $[0,\tau)$ where $\tau$ is a stopping time.  Let 
$(\langle M\rangle_t)_{0\leq t < \tau}$ be its quadratic variation and let $\langle M\rangle_{\tau} := \lim_{t\rightarrow \tau}\langle M\rangle_t \leq \infty$.
\begin{enumerate}
  \item
  On $\{ \langle M \rangle_{\tau}< \infty \}, \quad \lim_{t\rightarrow \tau} M_t $ exists a.s. in   $\R$.
  \item
  On $\{\langle M \rangle_{\tau}= \infty\},  \quad  \limsup_{t\rightarrow \tau}M_t = - \liminf_{t\rightarrow \tau }M_t = + \infty $ a.s.
 \end{enumerate}
\end{Lemme}
{\bf Proof}.
 (1)
  For any $p\geq 1$ let 
\[
  \sigma_p= \inf \{t\geq 0: \langle M \rangle_t\geq p\}. 
\]
In order every term in the following to be well defined we introduce a nondecreasing sequence of stopping times $(\tau_n)_{n\geq 1}$ with limit $\tau$ such that 
each stopped process  $(M_{\tau_n\wedge t})_{t\geq 0}$ is a uniformly integrable martingale for any $n\geq 1$.
The stopped process $(M_{t\wedge\tau_n\wedge\sigma_p})_{t\geq 0}$ is a $L^2$-bounded martingale and as $m,n \rightarrow \infty$
\[
  \E[(M_{\tau_n\wedge \sigma_p}- M_{\tau_m\wedge \sigma_p})^2]= \E[|\langle M \rangle_{\tau_n\wedge \sigma_p}-\langle M\rangle_{\tau_m\wedge\sigma_p}|] 
  \rightarrow 0.
\]
We set
  \[
  M^{(p)} := \lim_{n\rightarrow \infty} M_{\tau_n\wedge \sigma_p}= \lim_{t\rightarrow \tau} M_{t\wedge \sigma_p}.
\]
On $\{ \langle M \rangle_{\tau}<\infty\}$, the stopping times $\sigma_p$ are a.s. infinite from some $p$ on and we can set $M_{\tau}:=
\lim_{p\rightarrow \infty}M^{(p)}$.
Thus on this set
\[
  M_t\rightarrow M_{\tau} \quad { \rm as }  \quad t\rightarrow \tau .  
\]
(2) Let for any $t\geq 0$
\[
  T_t= \inf \{ 0\leq s\leq \tau : \langle M \rangle_s >t\}.
\]
There exist an enlargement $(\tilde{\Omega},(\tilde{\F}_t), \tilde{\P})$ of $(\Omega,(\F_{T_t}),\P) $ and a Brownian motion 
$\tilde{\beta}$ on $\tilde{\Omega}$ independent of $M$ such that the process
\[
B_t = \left\{
\begin{array}{lll}
      M_{T_t}  & \rm{if} & t< \langle M\rangle_{\tau} \\
      M_{\tau} + \tilde{\beta}_{t-\langle M \rangle_{\tau}} & \rm{if} & t\geq \langle M \rangle_{\tau}
\end{array}  \right .
\]
is a standard linear Brownian motion. As
\[
  \limsup_{t\rightarrow \infty} B_t = -\liminf_{t\rightarrow \infty} B_t = + \infty \quad \rm {a.s.}
\]
we obtain on $\{\langle M\rangle_{\tau}=\infty\}$
\[
  \limsup_{t\rightarrow \infty} M_{T_t} = - \liminf_{t\rightarrow \infty} M_{T_t} = + \infty
\]
and therefore 
\[  
\limsup_{t\rightarrow \tau} M_t = - \liminf_{t\rightarrow \tau}M_t = + \infty
 \].
 $\hfill \blacksquare$
 
 \vspace{0.4cm}

We will also  need the following consequence of the results in \cite{C95}  on multivalued
stochastic differential systems, completed with the  method used in \cite{CL97} and developed in \cite{L10} to check the lack of additional boundary process.
\begin{Prop}
\label{Th:inter}
Let $\alpha>0$, $\delta\geq 0$, ${\bf \sigma}=(\sigma_j^i; i,j=1,2) $ a $2\times2$-matrix, $(B,C)$ 
a Brownian motion in the plane, $b_1$ and $b_2$ two Lipschitz functions on $\R^2$,
$Z^1_0$ and $Z_0^2$ two ${\F}_0$-measurable  nonnegative random variables. There exists a unique strong  solution 
$(Z^1,Z^2)$ to the system
\begin{equation}
 \label{eq:mixte}
 \begin{array}{lll}
  Z^1_t & = & Z^1_0+ \sigma_1^1B_t+\sigma_2^1C_t + \alpha \int_0^t\frac{ds}{Z^1_s} + \int_0^tb_1(Z^1_s,Z^2_s)ds \\
  Z^2_t & = & Z^2_0+ \sigma_1^2B_t+\sigma_2^2C_t + \delta\int_0^t\frac{ds}{Z^2_s} + \int_0^tb_2(Z^1_s,Z^2_s)ds 
 \end{array}
\end{equation}
with the conditions $Z^1_t\geq 0$ if $\delta=0$ and $Z^1_t\geq 0,Z^2_t\geq 0$ if $\delta>0$.
\end{Prop}

\vspace{0.4cm}

It is worth noticing that the solutions to (\ref{eq:obli}) enjoy the Brownian scaling property. It means that if $(X,Y)$ is a solution to (\ref{eq:obli})
starting from $(X_0,Y_0)$ with driving Brownian motion $(B_t,C_t)$, then for any $c>0$ the process $(X_t^{\prime}:=c^{-1}X_{c^2t}, Y_t^{\prime}:=
c^{-1}Y_{c^2t}; t\geq 0)$ is a solution to (\ref{eq:obli}) starting from $(c^{-1}X_0,c^{-1}Y_0)$ with driving Brownian motion $(c^{-1}B_{c^2t},
c^{-1}C_{c^2t})$.

\vspace{0.5cm}

\section{Avoiding the corner}

We shall see in Section 5 that  existence and uniqueness of the solution to (\ref{eq:obli}) are easily obtained as soon 
as the solution process keeps away from the corner. Thus the question of attaining the corner in finite time is of great interest. 
For some class of reflection matrices, a necessary and sufficient condition is given in \cite{S15a}. 
Unfortunately, we are not able to provide such a complete answer and we have to be content with a collection of sufficient conditions 
ensuring nonattainability. We will just see at the end of Section 5 a degenerated case where  the corner is reached in finite time with full probability.

Our sufficient conditions are stated in the following  theorem.
In some cases (conditions $C_{2a}$ and $C_{2b}$ below), the comparison with Bessel processes suffices to conclude. In other
cases (conditions $C_1$ and $C_3$ below), we are first looking for a $C^2$-function $f$ on the punctured quadrant $S^{\bf 0}$ with limit $-\infty$
at the corner. Then we use  Lemma \ref{Th:kean} to show that $f(X_t,Y_t )$ cannot converge to $-\infty$ in finite time.

\begin{Th}
\label{Th:corn}
Let $(X,Y)$ be a solution to (\ref{eq:obli}). 
 We set 
\[
  \tau^{{\bf 0}}:= \inf \{t>0:(X_t,Y_t)={\bf 0}\}
\]
with the usual convention $\inf \emptyset =  \infty$.
Then $\P(\tau^{{\bf 0}}<\infty)=0$ if one of the following conditions is satisfied:
\begin{enumerate}
 \item
  $C_1: \beta \geq 0$, $\gamma \geq 0 $ and $-1\leq \rho\leq \alpha + \delta$.
 \item
  $C_{2a} : \alpha\geq \frac{1}{2}$ and  $\beta \geq 0$
 \item
  $C_{2b} : \delta\geq \frac{1}{2}$ and  $\gamma \geq 0$
 \item
 $C_3 :$ There exist $\lambda>0$ and $\mu>0$ such that
  \begin{itemize}
 \item
   $\lambda \alpha + \mu \gamma \geq 0$  
 \item
  $\lambda\beta+\mu \delta \geq 0$
   \item
   $(\sqrt{\lambda(\lambda\alpha+\mu \gamma)}+\sqrt{\mu(\lambda\beta +\mu \delta)})^2  \geq \frac{1}{2}(\lambda^2 + \mu^2+ 2 \rho \lambda \mu)$.
   \end{itemize}
 \end{enumerate}
\end{Th}
{\bf Proof}. \\
{\it Condition} $C_1$. For $\epsilon>0$ let 
\[
  \begin{array}{lll}
    \sigma^{\epsilon} & = & {\bf 1}_{\{(X_0,Y_0)={\bf 0}\}}\inf \{t>0: X_t+Y_t\geq \epsilon\} \\
    \tau^{{\bf 0},\epsilon} & = & \inf\{t>\sigma^{\epsilon}: (X_t,Y_t)={\bf 0}\} \,.
\end{array}
\]
As $\epsilon\downarrow 0$, $\sigma^{\epsilon}\downarrow 0$ 
and $ \tau^{{\bf 0},\epsilon}\downarrow \tau^{{\bf 0}}$.  We set $R_t= X_t^2+Y_t^2$. From It\^{o}'s formula we get 
for $t\in[\sigma^{\epsilon},\tau^{{\bf 0},\epsilon})$  
\[
  \begin{array}{ll}
   & \ln R_t \\
  = &  \ln R_{\sigma^{\epsilon}} + 2\int_{\sigma^{\epsilon}}^t\frac{X_s dB_s+Y_s dC_s}{R_s} + 2 \int_{\sigma^{\epsilon}}^t \frac{(\alpha+\delta)ds}{R_s}
     +2 \int_{\sigma^{\epsilon}}^t (\beta \frac{X_s}{Y_s}+\gamma \frac{Y_s}{X_s})\frac{ds}{R_s} -4 \rho \int_{\sigma^{\epsilon}}^t 
\frac{X_sY_s}{R_s^2}ds \\
 \geq &  \ln R_{\sigma^{\epsilon}} + M_t +  2\int_{\sigma^{\epsilon}}^t \frac{P_1(X_s,Y_s)}{R_s^2} ds
\end{array}
\]
where $M$ is a continuous local martingale on $[\sigma^{\epsilon},\tau^{{\bf 0},\epsilon})$ and $P_1(x,y)$ is the second degree homogeneous polynomial
\[
  P_1(x,y) = (\alpha+ \delta)(x^2+y^2) - 2\rho xy.
\]
Using Lemma \ref{Th:seconddegre} we check that this polynomial is nonnegative  on $S$  if $\rho \leq \alpha+ \delta$. Therefore
\[
  0\leq \int_{\sigma^{\epsilon}}^{\tau^{{\bf 0},\epsilon}} \frac{P_1(X_s,Y_s)}{R_s^2}ds \leq \infty
\]
 From Lemma \ref{Th:kean} we know
that depending on whether 
$\langle M\rangle_{\tau^{{\bf 0},\epsilon}}$ is finite or not, the local martingale $M_t$ either converges in $\R$ as $t\rightarrow \tau^{{\bf 0},\epsilon}$ 
or oscillates between $+\infty$ and $-\infty$. It cannot converge to $-\infty$. Thus $R_{\tau^{{\bf 0},\epsilon}}>0$ on 
 $\{\tau^{{\bf 0},\epsilon}<\infty\}$, which contradicts the definition of the moment $\tau^{{\bf 0},\epsilon}$. This proves that 
  $\tau^{{\bf 0},\epsilon}=\infty$ for every $\epsilon>0$, and therefore, $\tau^{{\bf 0}}=\infty$ a.s. \\
{\it Condition} $C_{2a}$ (resp. $C_{2b}$). Recall  the Bessel processes $U$ and $V$ in (\ref{eq:normal}). From Lemma 1 we get $X_t\geq U_t$ (resp. $Y_t\geq V_t$) and in this case $0$ is polar for 
$U$ (resp. V), so $X_t>0$ (resp. $Y_t>0$) for $t>0$. \\
{\it Condition} $C_3$.
We use again the notation for $\epsilon > 0$, $ \sigma^{\epsilon}$ and $ \tau^{{\bf 0},\epsilon}$.
We set $S_t=\lambda X_t+ \mu Y_t$ for $t\geq 0$, $\lambda>0$ and $\mu>0$. From It\^{o}'s formula we get
for $t\in[\sigma^{\epsilon},\tau^{{\bf 0},\epsilon})$
\[
 \begin{array}{ll}
  &\ln S_t \\
 = & \ln S_{\sigma^{\epsilon}} + \int_{\sigma^{\epsilon}}^t\frac{\lambda dB_s+\mu dC_s}{S_s} + 
  (\lambda\alpha+\mu \gamma)\int_{\sigma^{\epsilon}}^t\frac{ds}{X_sS_s}
  +(\lambda \beta+\mu \delta)\int_{\sigma^{\epsilon}}^t\frac{ds}{Y_sS_s}\\
  &  \qquad -\frac{1}{2}
     (\lambda^2 +\mu^2+ 2 \rho \lambda \mu)\int_{\sigma^{\epsilon}}^t\frac{ds}{S_s^2} \\
   = & \ln S_{\sigma^{\epsilon}} + M_t + \int_{\sigma^{\epsilon}}^t\frac{P_2(X_s,Y_s)}{X_sY_sS_s^2}ds
 \end{array}
\]
where $M$ is a continuous local martingale on  $[\sigma^{\epsilon}, \tau^{{\bf 0},\epsilon})$ and $P_2(x,y)$ is the second degree homogeneous polynomial
\[
  \begin{array}{ll}
 &   P_2(x,y)  \\
 = & \lambda(\lambda\beta+\mu\delta)x^2+\mu(\lambda\alpha+\mu\gamma)y^2+
  [\lambda(\lambda\alpha+\mu\gamma)+\mu(\lambda\beta+\mu\delta)-\frac{1}{2}(\lambda^2+\mu^2+ 2\rho \lambda \mu)]xy \,.
  \end{array}
\]
Using again Lemma \ref{Th:seconddegre}, we see that $P_2$ is nonnegative on $S$ if
\[
  \frac{1}{2}(\lambda^2+\mu^2+ 2\rho \lambda \mu)-  [\lambda(\lambda\alpha+\mu\gamma)+\mu(\lambda\beta+\mu\delta)]
  \leq 2 \sqrt{ \lambda(\lambda\beta+\mu\delta)\mu(\lambda\alpha+\mu\gamma)} .
\]
This is exactly  the  condition $C_3$. 
Therefore 
\[
  0\leq \int_{\sigma^{\epsilon}}^t \frac{P_2(X_s,Y_s)}{X_sY_sS_s^2}ds < \infty
\]
and so 
\[
0\leq \int_{\sigma^{\epsilon}}^{\tau^{{\bf 0},\epsilon}} \frac{P_2(X_s,Y_s)}{X_sY_sS_s^2}ds \leq \infty \,.
\]
Similarly, using   Lemma \ref{Th:kean} again, we see that the continuous local martingale $M$ either converges to a finite limit or oscillates 
between $+\infty$ and $-\infty$ when $t\rightarrow \tau^{{\bf 0},\epsilon}$.  It cannot converge to $-\infty$ and thus $S_{\tau^{{\bf 0},\epsilon}}>0$
on $\{\tau^{{\bf 0},\epsilon}<\infty\}$, proving $\P(\tau^{{\bf 0},\epsilon}<\infty)=0$. Letting finally  $\epsilon \rightarrow 0$ we obtain  
$\P(\tau^{{\bf 0}}<\infty)=0$.$\hfill  \blacksquare$ 

\vspace{0.5cm}
Condition $C_3$ is not explicit. We give two concrete examples when this condition holds true.

\begin{Cor} 
Assume $\rho=0$ , $\alpha=\delta$ and $|\beta|=|\gamma|$. Then the condition $C_3$ is satisfied  if
\begin{equation}
  \begin{array}{llll}
   \bullet & \beta^2\leq \alpha -\frac{1}{4} & \mbox{ when} &  \beta=-\gamma \\
   \bullet & -\beta \leq \alpha-\frac{1}{4} & \mbox{ when} & \beta=\gamma < 0\,.
  \end{array}
\end{equation}
\end{Cor}

{\bf Proof}. 
In both cases we take $\lambda=\mu>0$.\\
When $\beta=-\gamma\geq 0$,  the condition $C_3$  writes
 \begin{itemize}
  \item 
   $\alpha - \beta \geq 0 $
  \item
 $ 1 \leq (\sqrt{\alpha+\beta} + \sqrt{\alpha-\beta})^2 = 2 \alpha +2\sqrt{\alpha^2-\beta^2}$.
 \end{itemize}
If  $ \beta^2\leq \alpha -\frac{1}{4} $, then 
\[
  \alpha^2-\beta^2 \geq \alpha^2-\alpha+\frac{1}{4} \geq 0
\]
and
\[
  (1-2\alpha)^2=1-4\alpha+4\alpha^2 \leq 4(\alpha^2-\beta^2).
\]
When $\beta=\gamma<0$, the condition $C_3$  writes
 \begin{itemize}
  \item
   $\alpha+ \beta \geq 0$
  \item
   $1 \leq (2 \sqrt{\alpha+\beta})^2 = 4(\alpha + \beta)$
 \end{itemize}
and this is $-\beta \leq \alpha-\frac{1}{4}$. $\hfill  \blacksquare$

\vspace{0.3cm}

\begin{Cor}
\label{Th:asym}
Assume
\begin{itemize}
  \item 
  $\max\{\alpha,\delta\} \geq \frac{1}{2}$
  \item
  $2 \rho \leq \frac{\beta}{\delta}+ \frac{\gamma}{\alpha}$ \,.
\end{itemize}
Then  $\P(\tau^{{\bf 0}}<\infty)=0$ .
\end{Cor}

{\bf Proof}.
We may assume $\alpha\geq\frac{1}{2}$. If $\beta\geq 0$, then condition $C_{2a}$ holds true and the conclusion follows. If $\beta<0$, 
we will use condition $C_3$. We take $\lambda=\delta$ and $\mu=-\beta$. Then $\lambda \beta + \mu \delta=0$,
\[
  \begin{array}{lll}
  \lambda \alpha + \mu \gamma & = & \delta \alpha -\beta \gamma   \\
                                            &\geq & \delta \alpha -\beta \alpha (2 \rho- \frac{\beta}{\delta}) \\
                                        & = & \frac{\alpha}{\delta}(\delta^2+\beta^2-2\rho \beta \delta) \\
                                       & \geq & \frac{\alpha}{\delta}(\delta-\rho \beta)^2 \\
                                       & \geq & 0
   \end{array}
\]
and 
\[
  \begin{array}{ll}
   & (\sqrt{\lambda(\lambda\alpha+\mu \gamma)}+\sqrt{\mu(\lambda\beta +\mu \delta)})^2  - \frac{1}{2}(\lambda^2 + \mu^2 +2 \rho \lambda \mu) \\
  = & \lambda(\lambda\alpha+\mu \gamma) - \frac{1}{2}(\lambda^2 + \mu^2 +2 \rho \lambda \mu) \\
  \geq & \alpha (\delta^2+\beta^2-2\rho \beta \delta) -\frac{1}{2} (\delta^2+\beta^2-2\rho \beta \delta) \\
  \geq & (\alpha-\frac{1}{2}) (\delta-\rho \beta)^2 \\
  \geq &  0 \,.
  \end{array}
\]
$\hfill \blacksquare$

\vspace{0.5cm}

\section{Existence and uniqueness}

We now proceed to  the question of  existence and uniqueness of a global solution to (\ref{eq:obli}). We consider separately 
the three cases: first $\beta\geq 0$ and $\gamma \geq 0$, second $\beta \gamma <0$, third $\beta\leq 0$ and $\gamma\leq 0$. In the first and second cases, 
we construct the solution by switching from one edge to the other and patching the paths together. Thus it is essential to avoid the corner, as it was supposed 
 in \cite{S15b} in order to weakly approximate  an obliquely reflected Brownian motion. The third case uses a different proof and does not requires avoiding the corner. All three proofs  heavily use the comparison method of  Lemma \ref{Th:comp}.
 
\subsection{Case $\beta\geq 0$ and $\gamma \geq 0$}

\begin{Th}
 \label{Th:premier}
Assume $\beta\geq 0$, $\gamma \geq 0$ and one of the conditions $C_1$. $C_{2a}$, $C_{2b}$, $C_3$ is satisfied.
\begin{enumerate}
 \item
  There is a unique  solution to (\ref{eq:obli}) in $S^{{\bf 0}}$.
 \item
  There is a  solution to (\ref{eq:obli}) in $S$ starting from ${\bf 0}$.
 \item
  If $\alpha \delta\geq \beta \gamma$, there is a unique  solution  to (\ref{eq:obli}) in $S$.
\end{enumerate}
\end{Th}
{\bf Proof}. $1.$ Let $a>0$, $\epsilon>0$ and define for $(x,z)\in \R_+\times \R$
\[
  \psi_{\epsilon}(x,z):= \frac{1}{\max(\gamma x+z,\alpha \epsilon)}\,.
\]
This is a Lipschitz function. From Proposition \ref{Th:inter} we know that the system
\begin{equation}
  \label{eq:constru}
   \begin{array}{lll}
   X_t^{\epsilon} & = & X_0+B_t+\alpha\int_0^t\frac{ds}{X_s^{\epsilon}} +\alpha\beta\int_0^t \psi_{\epsilon}
  (X_s^{\epsilon},Z_s^{\epsilon})ds \geq 0 \\
  Z_t^{\epsilon}& = & -\gamma X_0+\alpha(Y_0+{\bf 1}_{\{Y_0=0\}}a)-\gamma B_t+\alpha C_t +\alpha(\alpha\delta-\beta\gamma)
    \int_0^t \psi_{\epsilon}(X_s^{\epsilon},Z_s^{\epsilon})ds
  \end{array}
\end{equation}
has a unique  solution. Let
\[
  \tau_Y^{\epsilon} := \inf \{t>0:\gamma X_t^{\epsilon}+ Z_t^{\epsilon}<\alpha \epsilon \} \,.
\]
If $0<\eta<\epsilon<a$ we deduce from the uniqueness that $(X^{\epsilon},Z^{\epsilon})$ and
$(X^{\eta},Z^{\eta})$ are identical on $ [0,\tau_Y^{\epsilon}]$. Patching together we can set
\[
 \begin{array}{lll}
  X_t & := & \lim_{\epsilon\rightarrow 0} X_t^{\epsilon} \\
  Y_t & := & \lim_{\epsilon\rightarrow 0} \frac{1}{\alpha}(\gamma X_t^{\epsilon} + Z_t^{\epsilon})
 \end{array}
\]
on 
$
\{(\omega,t)\in \Omega\times [0,\infty): Y_0(\omega)>0\; \rm{ and }\; 0\leq t<\tau_Y^0(\omega)\}, 
$
 where 
\[
  \tau_Y^0 := \lim_{\epsilon \rightarrow 0} \tau_Y^{\epsilon} \,.
\]
On this set, $(X,Y)$ is the unique solution to (\ref{eq:obli}). As we already noted, we have 
$X_t\geq U_t$ and $Y_t\geq V_t$. Therefore, on $\{Y_0>0\} \cap\{\tau_Y^0<\infty\}$,
\[
  \begin{array}{lll}
    \int_0^{\tau_Y^0} \frac{ds}{X_s} \leq \int_0^{\tau_Y^0} \frac{ds}{U_s} <\infty  &  \mbox{ and }  &
    \int_0^{\tau_Y^0} \frac{ds}{Y_s} \leq \int_0^{\tau_Y^0} \frac{ds}{V_s} <\infty 
  \end{array}
\]
and we can define
\begin{equation}
 \begin{array}{lllll}
  X_{\tau_Y^0} & := & \lim_{t\rightarrow \tau_Y^0} X_t & = & X_0+B_{\tau_Y^0} + \alpha \int_0^{\tau_Y^0}\frac{ds}{X_s}
       +\beta \int_0^{\tau_Y^0}\frac{ds}{Y_s} \\
  Y_{\tau_Y^0} & := & \lim_{t\rightarrow \tau_Y^0} Y_t & = & Y_0+C_{\tau_Y^0} + \gamma \int_0^{\tau_Y^0}\frac{ds}{X_s}
       +\delta \int_0^{\tau_Y^0}\frac{ds}{Y_s} \,.
  \end{array}
\end{equation}
We have  $Y_{\tau_Y^0}=0$. From Theorem \ref{Th:corn} we know that $X_{\tau_Y^0}>0$. 
In exactly the same way we can construct a solution 
on $\{Y_0>0\}$ in the interval $[T_1,T_2]$, where $T_1=\tau_Y^0$, $T_2= \inf \{t>T_1:X_t=0\}$. Iterating, we get a solution on 
$\{Y_0>0\} \times [0,\lim_{n\rightarrow\infty} T_n)$
where
\[
  \begin{array}{lll}
   T_{2p} & := & \inf\{t>T_{2p-1} : X_t=0\} \\
    T_{2p+1} & := & \inf\{t>T_{2p} : Y_t=0\}\,.
  \end{array}
\]
On $\{ Y_0>0\} \cap \{\lim_{n\rightarrow \infty}T_n<\infty\}$  
we  set $X_{\lim_{n\rightarrow \infty}T_n} := \lim_{p\rightarrow \infty} X_{T_{2p}}=0$ and 
$Y_{\lim_{n\rightarrow \infty}T_n} :=\lim_{p\rightarrow \infty} Y_{T_{2p+1}}=0$. The polarity of ${\bf 0}$ entails 
this is not possible in finite time and thus 
$\lim_{n\rightarrow \infty} T_n = \infty$. So we have obtained a unique global   solution on $\{Y_0>0\}$. In the same
way we obtain a unique global  solution on $\{X_0>0\}$ and as $\P((X_0,Y_0)= {\bf 0})=0$  the proof is complete.

2. Assume now $X_0=Y_0=0$. Let $(y_n)_{n\geq 1}$ be a sequence of real numbers (strictly) decreasing to $0$. From the above paragraph it follows there exists  
for any $n\geq 1$ a  unique solution $(X^n,Y^n)$ with values in $S^{\bf 0}$ to the system
\[
  \begin{array}{lll}
    X_t^n  & = & B_t+\alpha\int_0^t\frac{ds}{X_s^n} + \beta\int_0^t\frac{ds}{Y_s^n}  \\
    Y_t ^n& = & y_n+C_t+\gamma\int_0^t\frac{ds}{X^n_s} + \delta\int_0^t\frac{ds}{Y^n_s} \,.
 \end{array}
\]
 Let 
\[
  \tau:= \inf \{t>0: X^{n+1}_t<X_t^n\} \,.
\]
 Using  Lemma \ref{Th:comp}  we obtain $Y_t^{n+1} \leq Y_t^n$ on $[0,\tau]$. We note that
$(X_{\tau}^n,Y_{\tau}^n)\in S^{{\bf 0}}$ on $\{\tau<\infty\}$. On 
$\{Y_{\tau}^{n+1}=Y_{\tau}^n\}\cap \{\tau<\infty\}$, since $X_{\tau}^{n+1}=X_{\tau}^n$ and 
 the solution starting at time $\tau$ is unique, it follows 
that $X_t^{n+1}=X_t^n$ and $Y_t^{n+1}=Y_t^n$ on $[\tau, \infty)$.
On $\{Y_{\tau}^{n+1}<Y_{\tau}^n\}\cap \{\tau<\infty\}$, the continuity of solutions at time $\tau$ 
entails there exists $\rho>0$ such that $Y_t^{n+1}\leq Y_t^n$ on $[\tau,\tau+\rho]$. A second application of Lemma \ref{Th:comp}
 proves that $X_t^{n+1}\geq X_t^n$ on 
$[\tau,\tau+\rho]$, a contradiction
to the definition of $\tau$. Therefore $\P(\tau=\infty)=1$. It follows that 
$X_t^{n+1}\geq X_t^n$ and $Y_t^{n+1}\leq Y_t^n$ for any $t\in [0,\infty)$,
and we may define
\[
X_t:=\lim_{n\rightarrow\infty} \uparrow X_t^n  \qquad Y_t:=\lim_{n\rightarrow\infty} \downarrow Y_t^n \,.
\]
As $Y_t^n\geq V_t$ where $(U,V)$ is the solution to (\ref{eq:normal}) with $X_0=Y_0=0$, we have
\[
  \begin{array}{lll}
  X_t& = & B_t+ \alpha \lim_{n\rightarrow \infty} \int_0^t \frac{ds}{X_s^n}+ 
    \beta  \lim_{n\rightarrow \infty} \int_0^t \frac{ds}{Y_s^n} \\
  & =&  B_t + \alpha \int_0^t\frac{ds}{X_s}+ \beta \int_0^t\frac{ds}{Y_s} \\
   &< & \infty
  \end{array}
\]
and also 
\[
 \begin{array}{lll}
  Y_t& = &\lim_{n\rightarrow \infty}y_n +  C_t+ \gamma \lim_{n\rightarrow \infty} \int_0^t \frac{ds}{X_s^n}+ 
    \delta  \lim_{n\rightarrow \infty} \int_0^t \frac{ds}{Y_s^n} \\
  & =&  C_t + \gamma \int_0^t\frac{ds}{X_s}+ \delta \int_0^t\frac{ds}{Y_s} \\
   &< & \infty \,.
  \end{array}
\]

3. Assume finally $\alpha\delta-\beta \gamma \geq 0$. As the conclusion holds true if $\beta=\gamma=0$, we may also 
assume $\beta>0$. Let $(X,Y)$ be the solution to (\ref{eq:obli}) with $X_0=Y_0=0$ 
obtained in the previous paragraph and let $(X^{\prime},Y^{\prime})$ be another solution. Considering $(X^n,Y^n)$ again 
and replacing 
$(X^{n+1},Y^{n+1})$ with $(X^{\prime},Y^{\prime})$, the previous proof works and we finally obtain
$X^{\prime}_t\geq X_t$ and $Y^{\prime}_t\leq Y_t$. Then,
\begin{equation}
 \label{eq:uni} 
0 \leq  \delta(X^{\prime}_t-X_t)-\beta(Y_t^{\prime}-Y_t) = 
    \int_0^t (\alpha\delta-\beta\gamma)(\frac{1}{X_s^{\prime}}-\frac{1}{X_s})ds \leq 0 .
\end{equation}
Thus $X^{\prime}_t=X_t$ and $Y_t^{\prime}=Y_t$, proving  uniqueness. Replacing 
$(\delta,\beta)$ with $(\gamma,\alpha)$  in  equation 
(\ref{eq:uni}) we obtain the same conclusion 
if $\gamma>0$. $\hfill \blacksquare$

\vspace{0.3cm}

\begin{Rmq} {\rm  The statement in Theorem \ref{Th:premier} is not complete since
the problem of uniqueness when starting at the corner and $\alpha \delta<\beta \gamma$ is not solved. When considering the solution $(X,Y)$ 
in the above proof of existence, we have noted that  $X^{\prime}\geq X$ and $Y^{\prime}\leq Y$ for any other solution $(X^{\prime}, Y^{\prime})$. 
Thus uniqueness in law would be sufficient to obtain path uniqueness. A possible way to prove weak uniqueness could be the method in \cite{BP87}. This would be 
 far from our trajectorial methods and we don't go further in that direction. }
\end{Rmq}

\vspace{0.3cm}

\subsection{Case $\beta \gamma <0$}

\begin{Th}
\label{Th:deux}
Assume $\beta \gamma<0$ and one of the conditions $C_{2a}$ or $C_3$ is satisfied. Then, there
exists a unique  solution to (\ref{eq:obli}) in $S^{{\bf 0}}$.
\end{Th}
{\bf Proof}. 
Assume first  $\beta>0$, $\gamma<0$. The proof is similar to the proof of 1 in Theorem \ref{Th:premier}. The only change is that now 
$Y_t\leq V_t$. Therefore, on $\{Y_0>0\} \cap\{\tau_Y^0<\infty\}$,
\[
  \delta \int_0^{\tau_Y^0}\frac{ds}{Y_s}  \leq V_{\tau_Y^0} -Y_0 -C_{\tau_Y^0}-\gamma \int_0^{\tau_Y^0} \frac{ds}{U_s} < \infty
\]
and we can define $X_{\tau_Y^0}$ and  $Y_{\tau_Y^0}$ as previously done. The application of Theorem \ref{Th:corn} 
to the process on the time interval $[0, \tau_Y^0]$ shows that  $X_{\tau_Y^0}>0$ and we can iterate the construction as in Theorem 
\ref{Th:premier}. The proof if $\beta<0$, $\gamma>0$ is analogous.  $\hfill \blacksquare$

\vspace{0.3cm}

\subsection{Case $\beta\leq 0$ and $\gamma\leq 0$}

In this case we can give a full answer to the question of existence and uniqueness. When $|\rho|<1$, our condition of existence is 
analogous to  the condition found in \cite{W85} for the reflected Brownian in a wedge being a semimartingale, i.e. there is a convex combination 
of the directions of reflection that points into the wedge from the corner. It amounts to saying that the interaction matrix 
\[
  \left( \begin{array}{ll}
  \alpha & \beta \\
   \gamma & \delta 
   \end{array}
  \right)
\]
is completely-$\mathcal{S}$ in the terminology of  \cite{TW93, DW94, W95, DW96,  S15a, S15b, S15c}. 

\begin{Th}
\label{Th:trois}
Assume $\beta\leq 0$ and $\gamma\leq 0$.
\begin{enumerate}
  \item
   If $\alpha\delta>\beta\gamma$, there exists a unique  solution to (\ref{eq:obli}) in $S$.
  \item
  If $\alpha \delta \leq \beta \gamma$ and $ 1+\rho+|\alpha+\gamma|+|\beta+\delta|>0$, there is no solution.
  \item
  If $1+\rho=\alpha+\gamma=\beta+\delta=0$ and  $(X_0,Y_0)\neq{\bf 0}$ there exists a unique solution.
  \item
   If $1+\rho=\alpha+\gamma=\beta+\delta=0$ and  $(X_0,Y_0)={\bf 0}$ there is no solution.
\end{enumerate}
\end{Th}
{\bf Proof}.
1. Assume first $\alpha \delta > \beta \gamma$.\\
a) {\it Existence}.  Let $(h_n,n\geq 1)$ be a (strictly) increasing sequence of bounded positive 
nonincreasing Lipschitz functions converging to 
$1/x$ on $(0,\infty)$ and to $+\infty$ on $(-\infty,0]$. For instance we can take 
\[
  \begin{array}{lllll} 
      h_n(x) & = & (1-\frac{1}{n}) \frac{1}{x}  & \quad\mbox{ on } & [\frac{1}{n},\infty) \\
                 & = & n-1   & \quad \mbox{ on } & (-\infty,\frac{1}{n}] \,.
  \end{array}
\]
We consider for each $n\geq 1$ the system
\begin{equation}
 \label{eq:appro}
  \begin{array}{lll}
 X_t ^n& = & X_0+B_t+\alpha\int_0^t\frac{ds}{X_s^n} + \beta\int_0^th_n(Y_s^n)ds  \\
 Y_t^n & = & Y_0+C_t+\gamma\int_0^t h_n(X_s^n)ds + \delta\int_0^t\frac{ds}{Y_s^n} \,.
 \end{array}
\end{equation}
From Proposition \ref{Th:inter} it follows there exists a unique  solution to this system. 
 We set
\[
   \tau:= \inf\{s>0: X_s^{n+1}>X_s^n\}\,.
\]
We have $h_{n+1}(X_t^{n+1})\geq h_n(X_t^n)$ on $[0,\tau]$. A first application of Lemma \ref{Th:comp} shows that 
  $Y_t^{n+1}\leq Y_t^n$ on
$[0,\tau]$. Since
$h_{n+1}(Y_{\tau}^{n+1})> h_n(Y_{\tau}^n)$ on $\{\tau<\infty\}$, 
we deduce from the continuity of solutions that
there exists 
$\sigma>0$ such that $h_{n+1}(Y_t^{n+1})\geq h_n(Y_t^n)$ on $[\tau, \tau+\sigma]$. 
A second application of  Lemma \ref{Th:comp} shows that  $X_t^{n+1}\leq X_t^n$ on $[\tau, \tau+\sigma]$,
a contradiction to the definition of 
$\tau$. Thus $\P(\tau=\infty)=1$ proving that on the whole $[0,\infty)$ we have $X_t^{n+1}\leq X_t^n$ and $Y_t^{n+1}\leq Y_t^n$.
Then we can set for any $t\in [0,\infty)$
\[
    X_t:= \lim_{n\rightarrow \infty}X_t^n \quad \mbox{ and } \quad   Y_t:= \lim_{n\rightarrow \infty}Y_t^n \,.
\]
If $\alpha\delta>\beta \gamma$, there is a convex combination of the directions of repulsion pointing into 
the positive
quadrant, i.e. there exist $\lambda>0$ and $\mu>0$ such that $\lambda\alpha+\mu \gamma>0$ and $\mu\delta+\lambda
\beta>0$. For $n\geq 1$ and $t\geq 0$, 
\begin{equation}
 \label{eq:majo}
 \begin{array}{lll}
  \lambda U_t+\mu V_t & \geq & \lambda X_t^n+ \mu Y_t^n \\
    & \geq & \lambda X_0 + \mu Y_0 + \lambda B_t + \mu C_t + (\lambda\alpha+\mu\gamma) \int_0 ^t \frac{ds}{X_s^n}                    
                     + (\mu\delta+\lambda \beta) \int_0 ^t \frac{ds}{Y_s^n} \,.
  \end{array}
\end{equation}
Letting $n \rightarrow \infty$ in (\ref{eq:majo}) we obtain
\[
 \int_0^t \frac{ds}{X_s} <\infty \quad \mbox{and} \quad \int_0^t \frac{ds}{Y_s} <\infty \,.
\]
Then we may let $n$ go to $\infty$ in  (\ref{eq:appro}) proving that $(X,Y)$ is a solution to (\ref{eq:obli}).\\
b) {\it Uniqueness}. Let $(X^{\prime},Y^{\prime})$ be another solution to (\ref{eq:obli}). Replacing $(X^{n+1},Y^{n+1})$ with 
$(X^{\prime},Y^{\prime})$ we follow the above proof to obtain for $t\in[0,\infty)$ and $n\geq 1$ 
\[
   X_t^{\prime} \leq X^n_t \quad  \mbox{and} \quad  Y_t^{\prime} \leq Y^n_t
\]
Letting $n\rightarrow \infty$ we conclude
\[
  X_t^{\prime} \leq X_t \quad  \mbox{and} \quad  Y_t^{\prime} \leq Y_t \,.
\]
With the same $\lambda>0$ and $\mu>0$ as above,
\[
 0 \leq  \lambda (X_t-X_t^{\prime}) + \mu (Y_t-Y_t^{\prime})
  =  \int_0^t [(\lambda \alpha + \mu \gamma)
 (\frac{1}{X_s}-\frac{1}{X_s^{\prime}})+(\mu \delta+\lambda \beta)(\frac{1}{Y_s}-\frac{1}{Y_s^{\prime}})]ds \leq 0
\]
and therefore $X_t^{\prime}=X_t$, $Y_t^{\prime}=Y_t$. \\
2. If $\alpha \delta\leq \beta \gamma$ there exist $\lambda> 0$ and $\mu> 0$ such that $\lambda\alpha+\mu \gamma\leq 0$
and $\mu\delta+\lambda\beta\leq 0$. For that, just take
\begin{equation}
 \label{eq:encad}
 \frac{\alpha}{-\gamma}\leq \frac{\mu}{\lambda}\leq  \frac{-\beta}{\delta} .
\end{equation}
Let us consider the nonnegative quadratic $(\lambda^2+\mu^2+2\rho\lambda\mu)$. It is positive if $\rho\geq 0$. 
It is larger than $(1-\rho^2) \lambda^2>0$ if $-1<\rho<1$.  It may be positive
 if $\rho=-1$ and $|\alpha+\gamma|+|\beta+\delta|>0$ because from (\ref{eq:encad}) we may take $\lambda\neq \mu$.
 Thus, if $(X,Y)$ is a solution to (\ref{eq:obli}), for any $t\geq 0$,
 \[
   0 \leq  \lambda X_t + \mu Y_t \leq \lambda X_0 + \mu Y_0  + \lambda B_t + \mu C_t \,.
 \]
This is not possible since the paths of the  Brownian martingale 
 $(\lambda^2+\mu^2+2\rho \lambda\mu)^{-1/2} (\lambda B_t+\mu C_t)$ are not bounded below. 
 So there is no global solution.\\
3. Assume now $1+\rho=|\alpha+\gamma|+|\beta+\delta|=0$ and $X_0+Y_0>0$. 
The system becomes
\begin{equation}
 \label{eq:singu}
 \begin{array}{lllll}
 X_t & = & X_0+B_t+\alpha\int_0^t\frac{ds}{X_s} - \delta\int_0^t\frac{ds}{Y_s} & \geq & 0 \\
 Y_t & = & Y_0-B_t-\alpha\int_0^t\frac{ds}{X_s} + \delta\int_0^t\frac{ds}{Y_s} & \geq & 0.
 \end{array}
\end{equation}
 This entails for any $t\geq 0$
 \[
  X_t+Y_t=X_0+Y_0
\]
  and the first equation in (\ref{eq:singu}) reduces to
\begin{equation}
 \label{eq:simple}
 0 \leq  X_t  =  X_0+B_t+\alpha\int_0^t\frac{ds}{X_s} -\delta \int_0^t\frac{ds}{X_0+Y_0-X_s} \leq X_0+Y_0.  .
\end{equation}
Clearly pathwise uniqueness holds for equation (\ref{eq:simple})  since if there are two solutions $X$ and $X^{\prime}$, for any $t\geq 0$
\[
  \begin{array}{lll}
    (X_t-X_t^{\prime})^2 & =  & 2\alpha \int_0^t (X_s-X_s^{\prime})(\frac{1}{X_s}-\frac{1}{X_s^{\prime}})ds
                                          -2 \delta \int_0^t (X_s-X_s^{\prime})(\frac{1}{X_0+Y_0-X_s}-\frac{1}{X_0+Y_0-X_s^{\prime}})ds \\
                                         & \leq & 0 .
  \end{array}
\]
Consider now for any $c>0$ and $0\leq c_0 \leq c$ the equation 
\begin{equation}
\label{eq:intervalle}
   0 \leq Z_t  =  c_0+B_t+\alpha\int_0^t\frac{ds}{Z_s} - \delta\int_0^t\frac{ds}{c-Z_s}  \leq c .
\end{equation}
The solution $Z$ is a Brownian motion perturbed by a drift deriving from the concave potential   
\[
  \alpha \ln z  + \delta \ln (c-z)   \qquad {\rm for} \quad 0<z<c.
\]
From \cite{C95} we know that
equation (\ref{eq:intervalle}) has a unique solution living on the interval $[0,c]$. Thus there is a weak and then a unique strong solution to 
(\ref{eq:simple}). Setting 
\[
  Y_t=X_0+Y_0-X_t
\]
we obtain a unique strong solution to (\ref{eq:singu}).\\
4. We must have for any $t\geq 0$
\[
  X_t+Y_t=X_0+Y_0=0
\]
and thus $X_t=Y_t=0$.  But this is not possible for a solution to (\ref{eq:obli}). $\hfill \blacksquare$

\vspace{0.3cm}

We end this section by  considering a degenerate case where the solution hits the corner with probability one.

\begin{Prop}
\label{Th:collision}
Assume $\rho=1$, $\alpha\delta>\beta \gamma$, $\max\{\alpha,\delta\}<\frac{1}{2}$ and $\max\{\beta, \gamma\}\leq 0$. Then,
$\P(\tau^{{\bf 0}}<\infty)=1$.
\end{Prop}

{\bf Proof}.
We set
\[
 \begin{array}{lll}
  \tau_X^0 & := & \inf \{t>0: X_t=0\} \\
  \tau_Y^0 & := & \inf \{t>0: Y_t=0\} \,.
 \end{array}
\]
 The dimension of each Bessel process $U$ and $V$ in (\ref{eq:normal}) is less than 2, and $X\leq U$, $Y\leq V$.
 Then $\P( \tau_X^0<\infty=1)$ and $\P( \tau_Y^0<\infty=1)$.  Assume first $\alpha+\beta\leq \gamma+\delta$.
 Using  Lemma \ref{Th:comp} we obtain $X_t\leq Y_t \leq V_t$  on 
the time interval $[\tau_X^0, \infty)$. Therefore $\P(\tau^{{\bf 0}}<\infty)=1$. Same result when
$\alpha+\beta\geq \gamma+\delta$.  $\hfill  \blacksquare$

\vspace{0.4cm}

In the following pictures, we display the directions of interaction ${\bf  r_x}$ and ${\bf  r_y}$ defined in (\ref{eq:repul})
 in three illustrative instances.

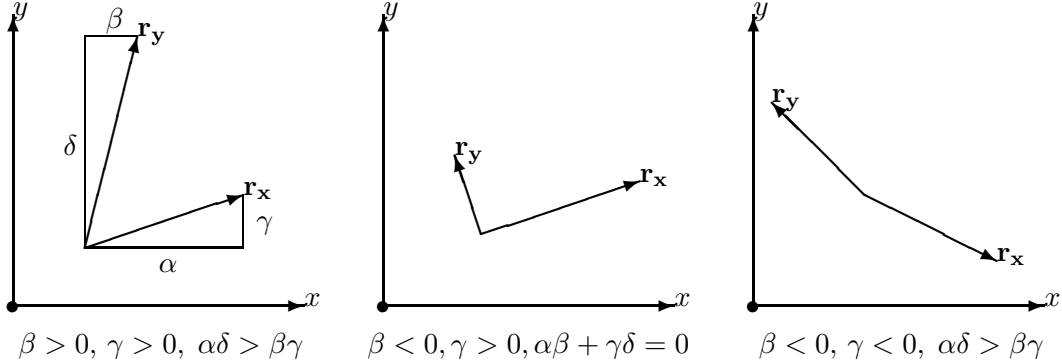
\begin{picture}(400,165)(5,20)
\begin{thicklines}
\put(5,40){$\bullet$}
\put(8,43){\vector(1,0){110}}
\put(8,43){\vector(0,1){110}}
\put(35,65){\vector(3,1){60}}
\put(35,65){\vector(1,4){20}}
\put(95,85){${\bf r_x}$}
\put(55,145){${\bf r_y}$}
\put(118,43){$x$}
\put(8,153){$y$}
\put(10,25){$\beta>0,\,\gamma>0,\;\alpha\delta> \beta\gamma$}
\put(145,40){$\bullet$}
\put(148,43){\vector(1,0){110}}
\put(148,43){\vector(0,1){110}}
\put(258,43){$x$}
\put(148,153){$y$}
\put(185,70){\vector(3,1){60}}
\put(185,70){\vector(-1,3){10}}
\put(245,90){${\bf r_x}$}
\put(175,100){${\bf r_y}$}
\put(142,25){$\beta<0,\gamma> 0,\alpha\beta+\gamma\delta=0$}
\put(285,40){$\bullet$}
\put(288,43){\vector(1,0){110}}
\put(288,43){\vector(0,1){110}}
\put(398,43){$x$}
\put(288,153){$y$}
\put(330,85){\vector(2,-1){50}}
\put(330,85){\vector(-1,1){35}}
\put(380,60){${\bf r_x}$}
\put(295,120){${\bf r_y}$}
\put(290,25){$\beta<0,\,\gamma<0,\;\alpha\delta>\beta\gamma$}

\end{thicklines}
\begin{thinlines}
\put(35,65){\line(0,1){80}}
\put(35,145){\line(1,0){20}}
\put(35,65){\line(1,0){60}}
\put(95,65){\line(0,1){20}}
\put(63,56){$\alpha$}
\put(100,72){$\gamma$}
\put(27,100){$\delta$}
\put(43,148){$\beta$}

\end{thinlines}
\end{picture}

\vspace{0.5cm}

\section{Avoiding the edges}

We now consider the question of  hitting  an  edge of the quadrant. Remember the definitions
\begin{equation}
 \begin{array}{lll}
  \tau_X^0 & := & \inf \{t>0: X_t=0\} \\
  \tau_Y^0 & := & \inf \{t>0: Y_t=0\} \,.
 \end{array}
\end{equation}
We already know that $\P(\tau_X^0<\infty)=0$ if $\alpha\geq\frac{1}{2}$ and $\beta\geq 0$.
Conversely, a comparison with the Bessel process $U$ shows 
that $\P(\tau_X^0<\infty)=1$ if $\alpha<\frac{1}{2}$ and $\beta\leq 0$. 
If we know that the 
corner is not hit and $\alpha\geq\frac{1}{2}$, we can get rid of the nonnegativity assumption on $\beta$. Since we are only interested in one coordinate,  we are looking for a function that is $C^2$ on $(0,\infty)$ and goes to $-\infty$ when approaching $0$. A natural candidate is the logarithmic function.

\begin{Prop}
\label{Th:side}
Assume $\P(\tau^{{\bf 0}}<\infty)=0$. If $\alpha\geq \frac{1}{2}$, then $\P(\tau^0_X<\infty)=0$.
\end{Prop}
{\bf Proof}. For $\eta>0$ let
\[
  \begin{array}{lll}
    \theta_X^{\eta} & = & {\bf 1}_{\{X_0=0\}}\inf \{t>0: X_t\geq \eta\} \\
    \tau_X^{0,\eta} & = & \inf\{t>\theta_X^{\eta}: X_t=0\} \,.
\end{array}
\]
As $\eta\downarrow 0$, $\theta^{\eta}_X\downarrow 0$  and $ \tau^{0,\eta}_X\downarrow \tau^0_X$.
For $t\in[\theta^{\eta}_X,\tau^{ 0,\eta}_X)$, from It\^{o}'s formula
\begin{equation}
\label{eq:bordX}
  \begin{array}{lll}
  \ln X_t & = & \ln X_{\theta_X^{\eta}} + \int_{\theta_X^{\eta}}^t \frac{dB_s}{X_s} +
    (\alpha-\frac{1}{2})\int_{\theta_X^{\eta}}^t \frac{ds}{X_s^2} + \beta \int_{\theta_X^{\eta}}^t \frac{ds}{X_sY_s}\,.
  \end{array}
\end{equation}
Since $\P(\tau^{{\bf 0}}<\infty)=0$,  on the set $\{\tau_X^{0,\eta}<\infty\}$ we have $Y_{\tau_X^{0,\eta}}>0$ because
$X_{\tau_X^{0,\eta}}=0$. On this set, from the definition of an O2BP, 
\[
  \int_{\theta_X^{\eta}}^{\tau_X^{0,\eta}} \frac{ds}{X_s} <\infty \;, 
 \qquad   \int_{\theta_X^{\eta}}^{\tau_X^{0,\eta}} \frac{ds}{Y_s} <\infty.
\]
As $Y_s>0$ on some interval $[\chi,\tau_X^{0,\eta}]$ with positive measure and $X_s>0$ on $[\theta_X^{\eta},\chi]$, we see that
\[
  | \beta |\int_{\theta_X^{\eta}}^{ \tau_X^{0,\eta}}\frac{ds}{X_sY_s}  < \infty \,.
\]
As $t\rightarrow \tau_X^{0,\eta}$, the local martingale in the r.h.s. of (\ref{eq:bordX}) cannot converge to $-\infty$. 
This entails that $\P(\tau_X^{0,\eta}<\infty)=0$ and therefore $\P(\tau^0_X<\infty)=0$. $\hfill \blacksquare$

\vspace{0.5cm}

We may also be interested in hitting either edge. This time we are looking for a function that is $C^2$  in the interior of $S$ and goes to 
$-\infty$ when approaching either edge. 

\begin{Prop}
\label{Th:sides}
Assume $\alpha\geq \frac{1}{2}$ and $\delta\geq \frac{1}{2}$. Then $\P(\tau^0_X<\infty)=\P(\tau^0_Y<\infty)=0$ if one of the following
conditions is satisfied:
 \begin{enumerate}
  \item
   $\beta\geq 0$
  \item
   $\gamma\geq 0$
  \item 
   $0 < \beta\gamma \leq (\alpha-\frac{1}{2})(\delta-\frac{1}{2})$.
 \end{enumerate}
\end{Prop}

{\bf Proof}. The proofs in the cases (1) and (2) ) are direct consequences of Theorem \ref{Th:corn} (under conditions $C_{2a}$ or $C_{2b}$)
and Proposition \ref{Th:side}. Assume now the conditions in (3) hold true. For $\epsilon>0$ let
 \[
   \begin{array}{lll}
    \sigma^{\epsilon} & = & {\bf 1}_{\{X_0Y_0=0\}} \inf\{t>0 : X_tY_t\geq \epsilon\} \\
    \tau^{\epsilon} & = & \inf\{t>\sigma^{\epsilon}:X_tY_t=0\}\;.
  \end{array}
\]
For $\lambda>0$ and $\mu>0$ we set
\[
  R_t=\lambda \ln X_t +\mu \ln Y_t.
\]
From It\^{o}'s formula we get for $t\in [\sigma^{\epsilon},\tau^{\epsilon})$
\[
  \begin{array}{lll}
   R_t & = & R_{\sigma^{\epsilon}} + \int_{\sigma^{\epsilon}}^t(\frac{\lambda}{X_s}dB_s+\frac{\mu}{Y_s}dC_s) + 
   \int_{\sigma^{\epsilon}}^t [\frac{\lambda (\alpha-\frac{1}{2})}
   {X_s^2} + \frac{\mu (\delta-\frac{1}{2})}{Y_s^2} + \frac{(\lambda\beta+\mu\gamma)}
   {X_sY_s} ]ds\\
    & = & R_{\sigma^{\epsilon}} + N_t +\int_{\sigma^{\epsilon}}^t \frac{P_3(X_s,Y_s)}{X_s^2Y_s^2}ds
  \end{array}
\]
where $N$ is a continuous local martingale and $P_3(x,y)$ is the second degree homogeneous polynomial 
\[
  P_3(x,y) = \mu(\delta-\frac{1}{2})x^2 + \lambda(\alpha-\frac{1}{2})y^2 +(\lambda \beta +\mu \gamma) xy.
\]
If the conditions (3) are satisfied,  we may take 
\[
  \begin{array}{lllll}
  \lambda & = & 2 (\alpha - \frac{1}{2})(\delta - \frac{1}{2}) - \beta \gamma & > &  0\\
  \mu  & = &  \beta^2 & > & 0. 
  \end{array}
\]
 Then, using again conditions (3), we  check that
\[
     P_3(x,y) 
  =  \beta^2 \left(\delta-\frac{1}{2}\right)x^2 + \left[2\left(\alpha-\frac{1}{2}\right)\left(\delta-\frac{1}{2}\right)-\beta \gamma\right]
\left(\alpha-\frac{1}{2}\right)y^2 
     + 2\left(\alpha-\frac{1}{2}\right)\left(\delta-\frac{1}{2}\right)\beta xy
\]
is  nonnegative on $S$. 
 The proof terminates as previously  in Theorem \ref{Th:corn} . 
$\hfill  \blacksquare$

\vspace{0.5cm}

\section{Associated local martingales}

We easily check that if $\max{\{\beta, \gamma,\alpha\delta-\beta\gamma\}}\geq 0$ and if $\rho>-1$, there exist some $\lambda\geq 0$
and $\mu\geq 0$ with  $\lambda  + \mu >0$ such that $\lambda X_t+\mu Y_t$ is not less than $\lambda (X_0 + B_t)+\mu (Y_0+C_t)$ which is 
 proportional to a real 
driftless Brownian motion. Therefore
\[
  \limsup_{t\rightarrow \infty} (\lambda X_t + \mu Y_t) = + \infty.
\]
An usual way in the study the asymptotic behavior  of the solutions to stochastic differential equations is introducing associated martingales. This is 
carried out through 
scale functions  (\cite{RY94} Section VII.3).
 There is no equivalent functions on the plane. However, in some particular 
cases we can find simple functions of O2BPs that are supermartingales or local martingales.

\begin{Prop}
Assume the following set of conditions:
\begin{equation}
\label{eq:super}
  \begin{array}{l}
   \alpha>1/2 \\
   \delta>1/2 \\
   \frac{\beta}{2\delta-1}+\frac{\gamma}{2\alpha-1} \geq  \rho  \\
   \rho>-1.  
   \end{array}
\end{equation}
Then 
\[
  M_t:=X_t^{1-2\alpha}Y_t^{1-2\delta}
\]
is  a positive supermartingale on $(0,\infty)$ which tends to $0$ as $t\rightarrow \infty$. It is a local martingale
if  the third inequality  in (\ref{eq:super})  is  an equality.  
\end{Prop}
{\bf Proof}. 
 For $\epsilon>0$ let 
\[
  \begin{array}{lll}
    \sigma^{\epsilon} & = & {\bf 1}_{\{X_0Y_0= 0\}}\inf \{t>0: X_tY_t\geq \epsilon\} \\
    \tau^{\epsilon} & = & \inf\{t>\sigma^{\epsilon}: X_tY_t= 0\} \,.
\end{array}
\]
As $\epsilon\downarrow 0$, $\sigma^{\epsilon}\downarrow 0$.
Up to $\tau^{\epsilon}$ , we get
from It\^{o}'s formula on $\{\sigma^{\epsilon}<\infty\}$
\[
  \begin{array}{lll}
  M_t & = &M_{\sigma^{\epsilon}} + \int_{\sigma^{\epsilon}}^t \frac{(1-2\alpha)Y_sdB_s+(1-2\delta)X_sdC_s}{X_s^{2\alpha} Y_s^{2\delta}} \\
        & & +   [ (1-2\alpha)\beta+(1-2\delta)\gamma +  \rho (2\alpha-1)(2\delta-1)]  \int_{\sigma^{\epsilon}}^t 
  \frac{ds}{X_s^{2\alpha} Y_s^{2 \delta}}.
\end{array}
\]
The function $f(x,y)=x^{1-2\alpha}y^{1-2\delta}$ is $C^2$ in the interior of $S$ and goes to $+\infty$ when approaching the edges of the quadrant,
whereas the finite variation part in the semimartingale decomposition of $M_t$ is nonincreasing.
We apply again Lemma \ref{Th:kean} and obtain that $\P(\tau^{\epsilon}<\infty)=0$. Then, letting $\epsilon \rightarrow 0$, we see that $M$ is a positive 
supermartingale on $(0,\infty)$. As such it tends to a limit $H\geq 0$ when $t\rightarrow \infty$ (\cite{RY94}, Corollary II.2.11 ). 
If $\beta\geq 0$ and $\gamma\geq 0$, then $X\geq U$ and $Y\geq V$ where $U$ and $V$ are the Bessel processes in  (\ref{eq:normal}). We have
 seen in Section 2 that $U_t\rightarrow \infty$ and $V_t \rightarrow \infty$ as $t\rightarrow \infty$; so $M_t\rightarrow 0$ as $t\rightarrow \infty$.
If now $\min\{\beta,\gamma\}<0$, say  $\gamma <0$, then we have 
  $Y\leq V$. 
From Corollary \ref{Th:besse} we deduce that for any $\epsilon>0$
\begin{equation}
 \label{eq:infini}
  \int_{\sigma^{\epsilon}}^{\infty}\frac{ds}{Y_s^2} \geq \int_{\sigma^{\epsilon}}^{\infty}\frac{ds}{V_s^2} =\infty .
\end{equation}
We
 consider the quadratic 
 variation $\langle M\rangle^{\epsilon}$ given by 
\[
 \begin{array}{lll}  
   \langle M\rangle_t^{\epsilon} &= &  \int_{\sigma^{\epsilon}}^t X_s^{-4\alpha} Y_s^{-4\delta}[(1-2\alpha)^2Y_s^2+(1-2\delta)^2X_s^2
    +2 \rho (1-2\alpha)(1-2\delta)X_sY_s]ds \\
             & = &  \int_{\sigma^{\epsilon}}^t M_s^2[\frac{(1-2\alpha)^2}{X_s^2}+\frac{(1-2\delta)^2}{Y_s^2}+2\rho\frac{(1-2\alpha)(1-2\delta)}{X_sY_s}]ds.
 \end{array}
\]
If $\rho\geq 0$,
\[
    \langle M\rangle_t^{\epsilon} \geq  \int_{\sigma^{\epsilon}}^t M_s^2\frac{(1-2\delta)^2}{Y_s^2} ds
\]
and if $\rho^2<1$,  
\[
   \langle M\rangle_t^{\epsilon} \geq (1-\rho^2) \int_{\sigma^{\epsilon}}^t M_s^2\frac{(1-2\delta)^2}{Y_s^2} ds.
\]
In both cases, using (\ref{eq:infini}) we should have $\langle M\rangle_{\infty}^{\epsilon}=\infty$ on the set $\{H>0\}\cap\{\sigma^{\epsilon}<\infty\}$, and then $\limsup_tM_t=-\liminf_tM_t=\infty$,
a contradiction with $M_t\rightarrow H$. Thus $H=0$ a.s. $\hfill \blacksquare$

\vspace{0.5cm}

During the proof we have seen that $\P(\tau_X^0<\infty)=\P(\tau_Y^0<\infty)=0$. In fact this is not a new result since here the conditions in Proposition 
\ref{Th:sides} are in force.
There is another case of interest.

\begin{Prop}
Assume $\alpha=\delta=1/2$, $\rho \beta\gamma<|\beta\gamma|$,  and there exists a solution to (\ref{eq:obli}) satisfying $\P(\tau^{{\bf 0}}<\infty)=0$.
Then
\[
 M_t:= \gamma \ln X_t -\beta \ln Y_t
\]
 is a continuous  local martingale on $(0,\infty)$ and  
 \[
   \limsup_{t\rightarrow \infty}{M_t}=-\liminf_{t\rightarrow\infty}{M_t}=\infty
\].
\end{Prop}
{\bf Proof}. We deduce  from Proposition \ref{Th:sides}  that $X_t>0$ and $ Y_t>0$ for any $t>0$. Let again for $\epsilon>0$ 
\[
 \sigma^{\epsilon}  =  {\bf 1}_{\{X_0Y_0= 0\}}\inf \{t>0: X_tY_t\geq \epsilon\}.
\]
Now It\^{o}'s formula gives
\[
  M_t=M_{\sigma^{\epsilon}} +\int_{\sigma^{\epsilon}}^t \left(\gamma \frac{dB_s}{X_s} -\beta \frac{dC_s}{Y_s}\right)
\]
and the quadratic variation is 
\[
  \langle M\rangle_t^{\epsilon} =   \int_{\sigma^{\epsilon}}^t\left(\frac{\gamma^2}{X_s^2} + 
    \frac{\beta^2}{Y_s^2}-2\rho\frac{\beta\gamma}{X_sY_s}\right) ds.
\]
Then, if $\rho\beta\gamma\leq 0$, 
 \[
    \langle M\rangle_t^{\epsilon} \geq  \int_{\sigma^{\epsilon}}^t\left(\frac{\gamma^2}{X_s^2} + 
    \frac{\beta^2}{Y_s^2}\right) ds,
\]
and if $\rho^2<1$,
\[
  \langle M\rangle_t^{\epsilon} \geq  (1-\rho^2) \max\left\{  \int_{\sigma^{\epsilon}}^t\frac{\gamma^2}{X_s^2}ds , 
     \int_{\sigma^{\epsilon}}^t\frac{\beta^2}{Y_s^2}ds \right\}.
\]
If $\beta>0$,  then $X\geq U$ where $U$ is the  Bessel process of dimension two in (\ref{eq:normal}) and 
 \[
   \limsup_{t\rightarrow \infty}{ (\ln X_t)} \geq  \limsup_{t\rightarrow \infty}{(\ln U_t)} = \infty
\]
If $\beta<0$, then $X\leq U$ and
 \[
   \liminf_{t\rightarrow \infty}{ (\ln X_t)} \leq  \liminf_{t\rightarrow \infty}{(\ln U_t)} =- \infty
\]
But
\[
  \ln X_t = \ln X_{\sigma^{\epsilon}} + \int_{\sigma^{\epsilon}}^t \frac{dB_s}{X_s} + \beta \int_{\sigma^{\epsilon}}^t \frac{ds}{X_sY_s}.
\]
and
\[ 
   \int_{\sigma^{\epsilon}}^t\frac{ds}{X_sY_s} \leq  \frac{1}{2} \int_{\sigma^{\epsilon}}^t \left( \frac{1}{X_s^2}+\frac{1}{Y_s^2}\right )ds.
\]
Thus  $\ln X_t$  would converge a.s. as $t\rightarrow \infty$ if
\[
 \int_{\sigma^{\epsilon}}^{\infty}\frac{ds}{X_s^2} < \infty  \quad \rm{and} \quad \int_{\sigma^{\epsilon}}^{\infty}\frac{ds}{Y_s^2} < \infty.
\]
Therefore $\langle M\rangle_{\infty}^{\epsilon}= \infty$ on $\{\sigma^{\epsilon}<\infty\}$ and the assertion is proved. $\hfill \blacksquare$
 
\vspace{0.5cm}

\section{Absolute continuity properties}

In this section we suppose $\rho=0$.
In some cases we easily obtain an absolute continuity property between the laws of O2BPs with various parameters. When  there exists a unique solution to
 (\ref{eq:obli}), we denote by $\P^{\alpha,\beta,\gamma,\delta}_{x,y}$ the  law on $C(\R_+,S)=C(\R_+,\R_+)\times C(\R_+,\R_+)$ 
of the solution starting at $(x,y)\in S$. We denote by $(U,V)$ the canonical map $(U_t(u,v),V_t(u,v))=(u(t),v(t))$ and by 
${\mathcal U}_t=\sigma((U_s,V_s),s\leq t)$ the canonical filtration. Recall that we denote by $\P^d_r$ the law of one-dimensional Bessel process of dimension $d$
starting at $r\geq0$.

\begin{Prop}
Assume $\rho=0$, $\delta\geq \frac{1}{2}$ and $y>0$. Then
\begin{equation}
 \P^{\alpha,\beta,0,\delta}_{x,y}\mid_{{\mathcal U}_t} = \exp{\left\{\beta\int_0^t\frac{dU_s}{V_s}-\alpha \beta \int_0^t\frac{ds}{U_sV_s}
   -\frac{\beta^2}{2} \int_0^t\frac{ds}{V_s^2}\right \}}\,.\,\P^{2\alpha+1}_x \otimes \P^{2\delta+1}_y \mid_{{\mathcal U}_t} .
\end{equation}
\end{Prop}
{\bf Proof}.
Under $\P^{2\alpha + 1}_x \otimes \P^{2\delta+1}_y$ the process 
\[
  B_t:= U_t-x- \alpha \int_0^t\frac{ds}{U_s}
\]
is a one-dimensional Brownian motion. The assumptions on $\delta$ and on $y$ imply that $V_t>0$ for any $t\geq 0$ and thus
\[
  \int_0^t\frac{ds}{v^2(s)} <\infty  \qquad \P_y^{2\delta+1}-a.s.
\]
Therefore, for $\P_y^{2\delta+1}$-almost every $v$, $\int_0^t\frac{dB_s}{v(s)}$ is a 
$\P_x^{2\alpha+1}$-centered Gaussian variable with variance $\int_0^t\frac{ds}{v^2(s)}$ and
\[
  \int\left(\int \exp{\left\{\beta\int_0^t\frac{dB_s}{v(s)}-\frac{\beta^2}{2}\int_0^t\frac{ds}{v(s)^2}\right\}}
   d\P^{2\alpha+1}_x(u) \right)d\P_y^{2\delta+1}(v) =1.
\]  
We see that
\[
 Z_t:=  \exp{\left\{\beta\int_0^t\frac{dB_s}{V_s}-\frac{\beta^2}{2}\int_0^t\frac{ds}{V_s^2}\right\}}
\]
is a $\P^{2\alpha+1}_x \otimes \P^{2\delta+1}_y$-positive martingale with expectation $1$. Setting  for any $T>0$
\[
   \Q_T := Z_T\, .\,\P^{2\alpha+1}_x \otimes \P^{2\delta+1}_y\mid_{{\mathcal U}_T}
\]
we check that under $\Q_T$ 
\[
  U_t-x-\alpha \int_0^t \frac{ds}{U_s}-\beta \int_0^t\frac{ds}{V_s}, \qquad 0\leq t\leq T
\]
is a real Brownian motion independent of the Brownian motion $\{V_t-y-\delta\int_0^t\frac{ds}{V_s},  \,0\leq t \leq T\}$.
Thus,
\[
  \Q_T= \P_{x,y}^{\alpha,\beta,0,\delta}\mid_{{\mathcal U}_T}. \qquad \qquad\qquad \qquad \hfill \blacksquare
\]

\vspace{0.4cm}

A second set of conditions is obtained using Novikov's criterion.

\begin{Prop}
Assume $\rho=0$ and the following set of conditions is satisfied:
\begin{equation}
 \begin{array}{ccc}
  |\beta|\leq \delta - 1/2 & & |\gamma|\leq \alpha- 1/2 \\
  x>0  & & y>0.
 \end{array}
\end{equation}
The process
\[
  Z_t:=  \exp{\left\{\beta\int_0^t\frac{dU_s}{V_s}- \gamma\int_0^t\frac{dV_s}{U_s}-(\alpha\beta+\gamma \delta) \int_0^t\frac{ds}{U_sV_s}
-\frac{\beta^2}{2}\int_0^t\frac{ds}{V_s^2}-\frac{\gamma^2}{2}\int_0^t\frac{ds}{U_s^2}\right\}},
\]
 is a $\P^{2\alpha+1}_x \otimes \P^{2\delta+1}_y$-positive martingale with expectation $1$ and
\begin{equation}
  \P^{\alpha,\beta,\gamma,\delta}_{x,y}\mid_{{\mathcal U}_t} = Z_t \,. \, \P^{2\alpha+1}_x \otimes \P^{2\delta+1}_y \mid_{{\mathcal U}_t}.
\end{equation}
\end{Prop}
{\bf Proof}.
Under $ \P^{2\alpha+1}_x \otimes \P^{2\delta+1}_y$ the processes
\[
  B_t:= U_t-x- \alpha \int_0^t\frac{ds}{U_s}
\]  
  and
\[
  C_t:= V_t-y- \delta \int_0^t\frac{ds}{V_s}
\]
are independent Brownian motions. Then $Z_t$ may be written
\[
  Z_t=   \exp{\left\{\beta\int_0^t\frac{dB_s}{V_s}+\gamma\int_0^t\frac{dC_s}{U_s}
  -\frac{\beta^2}{2}\int_0^t\frac{ds}{V_s^2}-\frac{\gamma^2}{2}\int_0^t\frac{ds}{U_s^2}\right\}},
\]
which shows that $Z_t$ is a positive local martingale. Using (\ref{eq:absocon}) we compute
\[
  \begin{array}{lll}
  \int  \int \exp{\{\frac{\gamma^2}{2}\int_0^t\frac{ds}{U_s^2}\}}\; d\P^{2\alpha+1}_x \otimes  d\P^{2\delta+1}_y & = & 
    \int  \exp{\{\frac{\gamma^2}{2}\int_0^t\frac{ds}{R_s^2}\}}\; d\P^{2\alpha+1}_x \\
     & = & \int (\frac{R_t}{x})^{\alpha-1/2} \exp {\{-\frac{(\alpha-1/2)^2-\gamma^2}{2}\int_0^t\frac{ds}{R_s^2} \}}\; d\P^2_x \\
    & <& \infty
   \end{array}
\]
since a Bessel process of dimension two has finite moments of any order. Finally,
\[
  \begin{array}{ll}
   & \int \int \exp{\{\frac{\beta^2}{2}\int_0^t\frac{ds}{V_s^2}\}}\exp{\{\frac{\gamma^2}{2}\int_0^t\frac{ds}{U_s^2}\}}\; d\P^{2\alpha+1}_x \otimes d\P^{2\delta+1}_y \\
   = & ( \int  \exp{\{\frac{\beta^2}{2}\int_0^t\frac{ds}{R_s^2}\}}\; d\P^{2\delta+1}_x)
  ( \int  \exp{\{\frac{\gamma^2}{2}\int_0^t\frac{ds}{R_s^2}\}}\; d\P^{2\alpha+1}_x)\\
    < & \infty
  \end{array}
\]
and Novikov's criterion (\cite{RY94}, Proposition VIII.1.15) proves that $Z_t$ has expectation $1$ with respect to
$\P^{2\alpha+1}_x \otimes  \P^{2\delta+1}_y $ . We easily see  that
$\{B_t-\beta\int_0^t \frac{ds}{V_s}, 0\leq t \leq T\}$ and
$\{C_t-\gamma\int_0^t \frac{ds}{U_s}, 0\leq t \leq T\}$ are independent Brownian motions under the probability with density $Z_T$ and 
this proves that the new probability is  $\P^{\alpha,\beta,\gamma,\delta}_{x,y}\mid_{{\mathcal U}_T}$. $\hfill \blacksquare$

\vspace{0.5cm}

\section{Product form stationary distribution}

 For $a>0$ and $c>0$, let  $\Gamma(a,c)$  be the probability measure on 
$[0,\infty)$ with density 
\[
\gamma(x;a,c):= \frac{c^a}{\Gamma(a)} x^{a-1}e^{-cx} 
\]
and characteristic function
\begin{equation}
\label{eq:caract}
  \phi(\lambda;a,c) := \int_0^{\infty} e^{i\lambda x} \gamma(x;a,c) dx = \left(1-\frac{i\lambda}{c}\right)^{-a}\,.
\end{equation}
For $a>0$, $b>0$, let  $B(a,b)$  be the probability measure   on $[0,1]$ with density
\[
  \beta(x;a,b) := \frac{\Gamma(a+b)}{\Gamma(a)\Gamma(b)} x^{a-1}(1-x)^{b-1}.
\]

It is well-known that the function $g(x)=x^{d-1}$ is an invariant measure density 
on $[0, \infty)$ for the Bessel process of dimension $d$. To get a stationary probability 
we have to 
introduce a negative drift that entails  positive recurrence. Then, the process
\[
  X_t=X_0 +B_t + \frac{d-1}{2}\int_0^t\frac{ds}{X_s} - \theta t,
\]
with $\theta>0$, has $\gamma(x;d,2 \theta)$ as stationary density. As for the  Brownian motion reflected at $0$ with constant drift $-\theta$,  it has the stationary exponential density $2\theta \exp\{-2 \theta x\}$. 
In the bidimensional case, the drifted obliquely reflected Brownian motion  has a stationary density in the form of product of two exponential densities if and only if it satisfies at once \cite{HW87a, W95}:
\begin{itemize}
 \item
  an invertibility condition on the reflection matrix;
 \item
  a positivity condition on the exponential coefficients;
 \item
  a skew-symmetry condition.
\end{itemize}
Therefore it is natural to ask wether one can find similar conditions for drifted O2BPs ensuring existence of a stationary distribution in the form of product of two 
gamma distributions. The answer is positive and given in the next theorem. In fact, this may be considered as a consequence of the study in \cite{OO14} which introduced a generalised reflected Brownian motion associated with a potential $U$ regular on the whole real line. The difference is that here the logarithmic potential
 is defined only on the positive axis. The proof below is an adaptation of the proof in \cite{OO14} to this special case.

We introduce an additional constant drift $(-\theta,-\eta)$ in order to make the solution a  recurrent process  in the nonnegative quadrant. We consider the system
\begin{equation}
 \label{eq:derive}
 \begin{array}{lll}
 X_t & = & X_0 + B_t + \alpha\int_0^t\frac{ds}{X_s} + \beta \int_0^t\frac{ds}{Y_s} - \theta t\\
 Y_t & = & Y_0 + C_t + \gamma\int_0^t\frac{ds}{X_s} + \delta \int_0^t\frac{ds}{Y_s} - \eta t
 \end{array}
\end{equation}
with the conditions $X_t\geq 0$, $Y_t\geq 0$. Changing probability through a Girsanov transformation, we easily check that Theorem \ref{Th:corn} in Section 4
 is still valid for this drifted system. Moreover, the proofs  in Section 5 do not bother whether the Brownian motions $B$ and $C$ 
are drifted or not drifted. Therefore
existence and uniqueness results in Section  5  hold true for the solution to (\ref{eq:derive}). 

\begin{Th}
\label{Th:produit}
Assume there exists a unique solution  to (\ref{eq:derive}) in $S^{{\bf 0}}$ or $S$. This process has 
an invariant distribution in the form $\Gamma(a,c)\otimes \Gamma(b,d)$ with $a>1$ and $b>1$  if and only if at once
\[
 \begin{array}{lllll}
  \bullet & \alpha\delta-\beta\gamma \neq 0 & & &  \mbox{(invertibility of the interaction matrix})\\
  \bullet &\alpha\eta-\gamma\theta>0 & \mbox{ and } &  \delta \theta-\beta \eta >0 & \mbox{(positivity of exponents)} \\
         \bullet & 2\rho=\frac{\beta}{\delta}+\frac{\gamma}{\alpha} & & & \mbox{(skew-symmetry)}
 \end{array}
\]
Under these conditions, the unique solution is given by
\[
  \begin{array}{llll}
  \bullet & a & = & 1 + 2 \alpha \\
  \bullet & b & = & 1 + 2 \delta \\
  \bullet & c & = & 2\alpha\frac{\delta \theta-\beta \eta}{\alpha\delta-\beta\gamma} \\
  \bullet & d & = & 2 \delta  \frac{\alpha\eta-\gamma\theta}{\alpha\delta-\beta\gamma}.
 \end{array}
\]
\end{Th} 
 {\bf Proof}. We first remark that under the above conditions $\alpha \delta - \beta \gamma$ cannot be $<0$, because in that case there is no solution if 
$\beta$ and $\gamma$ are $<0 $ and if $\beta$ and $\gamma$ are $>0$ the skew-symmetry condition is not satisfied since 
$2 \rho \leq 2< \beta/\delta + \gamma/\alpha$. So we could have written 
$\alpha \delta -\beta \gamma >0$ as first condition.
 Let now
\[
  q(x,y)= \gamma(x;a,c) \gamma(y;b,d) \qquad \mbox{for } x\geq 0,y\geq 0 \,.
\]
 The infinitesimal generator of the diffusion (\ref{eq:derive}) is given by
\[
  L= \frac{1}{2}(\frac{\partial^2}{\partial x^2}+\frac{\partial^2}{\partial y^2}+2\rho \frac{\partial^2}{\partial x \partial y})
    +(\frac{\alpha}{x}+\frac{\beta}{y}-\theta)\frac{\partial}{\partial x}
    +(\frac{\gamma}{x}+\frac{\delta}{y}-\eta)\frac{\partial}{\partial y} \,.
\]
Assume $(X_0,Y_0)$ has density $q$ and characteristic function
\[
   \E[e^{i(\lambda X_0 +\mu Y_0)}] = \phi(\lambda;a,c)\, \phi(\mu;b,d) \;.
\]
We set $f(x,y)=e^{i(\lambda x + \mu y)}$. We want to prove that for any $t\geq 0$,  $\lambda$, and $\mu$,
\[
  \E[f(X_t,Y_t)] = \E[f(X_0,Y_0)]\;( =  \phi(\lambda;a,c) \phi(\mu;b,d) \;) \;.
\]
It is enough to prove that 
\[
  \int_0^{\infty}\int_0^{\infty}Lf(x,y)q(x,y)\,dxdy\,= 0 
\]
for any $\lambda$ and $\mu$. Let 
\[
  R(\lambda, \mu) :=  [-\frac{1}{2}(\lambda^2+\mu^2 +2\rho \lambda \mu)+i\lambda(\frac{\alpha}{x}+\frac{\beta}{y}-\theta) +
      i\mu(\frac{\gamma}{x}+\frac{\delta}{y}-\eta)]
\]
and  compute
\[ 
  \begin{array}{lll}    
 \int_0^{\infty}\int_0^{\infty}Lf(x,y)q(x,y)\,dxdy &
 =  & \int_0^{\infty}\int_0^{\infty} R(\lambda, \mu)    e^{i(\lambda x+\mu y)}  \gamma(x;a,c)\gamma(y;b,d)\;dxdy \\
  & = & S(\lambda, \mu) \,  \phi(\lambda;a,c) \, \phi(\mu;b,d) 
  \end{array}
\]
where $S(\lambda,\mu)$ is a second degree polynomial. In the last computation we used formula (\ref{eq:caract}) several times.
Setting to zero the coefficients  of  $\lambda$, $\mu$, $\lambda^2$, $\mu^2$ and 
$\lambda \mu$ in polynomial $S$ 
we obtain the set of conditions:
\[ 
 \begin{array}{lll}
   0 & = & -\theta+\frac{\alpha c}{a-1}+ \frac{\beta d}{b-1} \\
  0  & = & -\eta + \frac{\gamma c}{a-1} + \frac{\delta d}{b-1} \\
  0 & = & -\frac{1}{2}+ \frac{\alpha}{a-1} \\
  0 & = & -\frac{1}{2}+ \frac{\delta}{b-1} \\
  0 & = & -\rho + \frac{\gamma}{a-1}+ \frac{\beta}{b-1} \;.
 \end{array}
\]

The solution is given by the specified values for $a,b,c,d$ and the skew-symmetry condition.  $\hfill \blacksquare$ 
\vspace{0.4cm}

It was  proved in \cite{W87} (see also \cite{S15a}) that under a skew-symmetry condition the obliquely reflected Brownian motion does not reach the 
non-smooth part of the boundary. In the same way the above 
  skew-symmetry equality is reminiscent of the second condition in Corollary \ref{Th:asym}, which is an inequality.
 We therefore get a partial but  handy statement.

\begin{Cor}
\label{Th:prod}
Assume   
\[
 \begin{array}{llll}
  \bullet & \max\{\alpha,\delta\} \geq \frac{1}{2} &&  \\
  \bullet & \alpha\delta-\beta\gamma > 0 & &   \\
  \bullet & \alpha\eta-\gamma\theta>0 & \mbox{ and } & \delta \theta-\beta \eta >0
      \\
  \bullet & 2\rho=\frac{\beta}{\delta}+\frac{\gamma}{\alpha}. & & 
 \end{array}
\]
Then the unique solution to (\ref{eq:derive})  has an invariant distribution given by 
\[
  \Gamma\left(1+2\alpha, 2\alpha\frac{\delta \theta-\beta \eta}{\alpha\delta-\beta\gamma}\right) \otimes
 \Gamma \left(1+2\delta, 2 \delta  \frac{\alpha\eta-\gamma\theta}{\alpha\delta-\beta\gamma}\right).
\]
\end{Cor}

{\bf Proof}. This is a straightforward consequence of Corollary  \ref{Th:asym}, Theorems \ref{Th:premier},
\ref{Th:deux}, \ref{Th:trois} and \ref{Th:produit}. $\hfill \blacksquare$

\vspace{0.3cm}

\begin{Rmq}
{\rm  When $\rho=0$, the skew-symmetry condition  $\alpha \beta + \gamma \delta =0$ means that the directions of interaction 
${\bf r_x}$ and ${\bf r_y}$  are orthogonal. Then the second condition in Theorem  \ref{Th:produit}  means that $(\theta,\eta)$ points into the 
 interior of the quadrant designed by ${\bf r_x}$ and  ${\bf r_y}$. This condition ensures recurrence of the process, while the first condition is now a consequence of the 
  skew-symmetry condition. }
\end{Rmq}

\begin{Rmq}
{\rm With the same proof, we may check that under the skew-symmetry condition, 
when $\theta=\eta=0$, the function $q(x,y)=x^{2\alpha}y^{2\delta}$ is a non-integrable invariant density that  does not depend 
on the other parameters.}
\end{Rmq}

\vspace{0.4cm}

A simple change of variables (beta-gamma algebra) provides  the following result.

\begin{Cor}
With the conditions and notations of Theorem \ref{Th:produit}, the two-dimensional process
\begin{equation}
  \begin{array}{lll}
    W_t & := & \frac{cX_t}{cX_t+dY_t} \\
    Z_t & := & cX_t+dY_t
  \end{array}
\end{equation}
has  $B(a,b)\otimes \Gamma(a+b,1)$ for invariant distribution on $[0,1]\times [0,\infty)$.
\end{Cor}

\vspace{0.4cm}

{\bf Acknowledgments.} The author is deeply indebted to the anonymous referee for careful reading, encouragements and sound advice. They have given rise to
 fundamental improvements in the content of this paper.

\end{document}